\documentclass[a4paper,12pt]{amsart}

\usepackage{amsmath,amssymb,amsthm,graphicx,pstricks,xypic,comment,fancyhdr}
\usepackage[all]{xy}
\usepackage{epic}
\usepackage{eufrak}
\usepackage{lscape}

\newcommand{\Hom}{{\sf Hom }}
\newcommand{\End}{{\sf End }}
\newcommand{\Ext}{{\sf Ext }}

\newcommand{\Coker}{{\sf Coker }}

\renewcommand{\mod}{{\sf mod \hspace{.02in}  }}
\newcommand{\Mod}{{\sf Mod \hspace{.02in} }}

\newcommand{\per}{{\sf per \hspace{.02in}  }}

\newcommand{\Sub}{{\sf Sub \hspace{.02in} }}

\newcommand{\add}{{\sf add \hspace{.02in} }}
\newcommand{\fl}{{\sf f.l. \hspace{.02in} }}
\newcommand{\gr}{\sf gr \hspace{.02in} }
\newcommand{\Jac}{{\sf Jac} }

\newcommand{\ten}{\otimes}
\newcommand{\lten}{\overset{\mathbf{L}}{\ten}}
\newcommand{\RHom}{\mathbf{R}{\sf Hom }}

\newcommand{\A}{\bar{A}}\newcommand{\B}{\bar{B}}

\newcommand{\Cc}{\mathcal{C}}
\newcommand{\Dd}{\mathcal{D}}
\newcommand{\Ee}{\mathcal{E}}

\newcommand{\Ii}{\mathcal{I}}

\newcommand{\Pp}{\mathcal{P}}

\newcommand{\ww}{\mathbf{w}}
\newcommand{\BPi}{{\mathbf\Pi}}

\newcommand{\bsm}{\begin{smallmatrix}}
\newcommand{\esm}{\end{smallmatrix}}

\newtheorem{thma}{Theorem}[section]
\newtheorem*{thm*}{Th{\'e}or{\`e}me}
\newtheorem{lema}[thma]{Lemma}
\newtheorem*{lem*}{Lemme}

\newtheorem{cora}[thma]{Corollary}
\newtheorem*{prop*}{Proposition}
\newtheorem{prop}[thma]{Proposition}
\theoremstyle{remark}

\newtheorem{rema}[thma]{Remark}

\theoremstyle{definition}

\newtheorem{dfa}[thma]{Definition}

\title{The ubiquity of generalized cluster categories}
\author{Claire Amiot}
\address{Institut de Recherche Math\'ematique Avanc\'ee, 7 rue Ren\'e Descartes, 67000 Strasbourg, France}
\email{amiot@math.unistra.fr}
\thanks{All authors were supported by the Storforsk-grant 167130 from the Norwegian Research Council}
\author{Idun Reiten}
\address{Insitutt for matematiske fag,
Norges Teknisk-Naturvitenskapelige Universitet,
N-7491 Trondheim, Norway}
\email{idun.reiten@math.ntnu.no}
\author{Gordana Todorov}
\address{Department of Mathematics, Northeastern University, Boston, MA 02115, USA}
\email{g.todorov@neu.edu}\thanks{The third author was also supported by the NSA-grant MSPF-08G-228}
\begin{document}

\maketitle

\bigskip
\begin{abstract}
Associated with a finite dimensional algebra of global dimension at most 2, a generalized cluster category was introduced in \cite{Ami3}. It was shown to be triangulated, and 2-Calabi-Yau when it is $\Hom$-finite. By definition, the cluster categories of \cite{Bua} are a special case. In this paper we show that a large class of 2-Calabi-Yau triangulated categories, including those associated with elements in Coxeter groups from \cite{Bua2}, are triangle equivalent to generalized cluster categories. This was already shown for some special elements in \cite{Ami3} and then more generally for $c$-sortable elements in~\cite{AIRT}.
\end{abstract}

\tableofcontents
\section*{Introduction}
Throughout this paper $k$ is an algebraically closed field.
Let $Q$ be a finite quiver without oriented cycles. In \cite{Bua}, the cluster category $\Cc_Q$ was defined
to be the orbit category $\Dd^{\rm b}(kQ)/\tau^-[1]$, where $\tau$ is the
Auslander-Reiten translation in the bounded derived category
$\Dd^{\rm b}(kQ)$. The category
$\Cc_Q$ is $\Hom$-finite, triangulated \cite{Kel}, and
$2$-Calabi-Yau (2-CY for short), that is, there is a functorial isomorphism
$D\Hom_{\Cc_Q}(X,Y)\simeq \Hom_{\Cc_Q}(Y,X[2])$, where $D=\Hom_k(-,k)$.
A theory for a special kind of objects, called \emph{cluster-tilting objects}, was developed in \cite{Bua}.
This work was motivated via \cite{Mar} by the Fomin-Zelevinsky
theory of cluster algebras \cite{FZ1}, where the cluster-tilting objects are the analogs of clusters.
  
Another category where a similar theory was developed is the category $\mod \Lambda$ 
of finite dimensional modules over a preprojective algebra $\Lambda$ of Dynkin type \cite{Gei3,Gei}. 
This category is $\Hom$-finite and Frobenius. 
Moreover, it is stably 2-CY, that is, its stable category $\underline{\mod}\Lambda$ (which is triangulated) is  $2$-CY.  

Much of the work on cluster categories from \cite{Bua, BMR07, BMR08} has been generalized to the setting of 2-CY triangulated
 categories with cluster-tilting objects, and new results have been proved in the general setting ( \cite{Iya,Kel4,Kel5}, and others).
It is of interest to investigate such categories, both for developing new theory and for providing 
applications to new classes of cluster algebras. In particular, it is of interest to find classes 
of 2-CY triangulated categories with cluster-tilting objects. An important class is the stable categories 
$\underline{\Ee}_w$ of the Frobenius categories $\Ee_w$ associated with elements $w$ in Coxeter groups 
\cite{Bua2}, (see  \cite{Gei4} for independent work when $w$ is adaptable). This class contains both the cluster categories $\Cc_Q$ and $\underline{\mod}\Lambda$  discussed above as special cases 
(see \cite{Bua2}, \cite{Gei4}). In $\Ee_w$ and $\underline{\Ee}_w$, there are \emph{standard} cluster-tilting objects $T_{\ww}$ associated with any reduced expression $\ww$ of $w$.

A new class of triangulated 2-CY categories was introduced in \cite{Ami3}.
They are \emph{generalized cluster categories} $\Cc_{\A}$ associated with algebras $\A$
of global dimension at most 2, rather than global dimension 1. In this case the orbit category 
$\Dd^{\rm b}(\A)/\tau^-[1]$ is not necessarily triangulated, so $\Cc_{\A}$ is defined to be its triangulated hull.
If $\Cc_{\A}$ is $\Hom$-finite, then it is triangulated 2-CY and $\A$ is a cluster-tilting object in $\Cc_{\A}$.

A natural question is how the generalized cluster categories are related to the previous classes of $\Hom$-finite triangulated 2-CY categories. 
It was already shown in \cite{Ami3} that some classes of categories $\underline{\Ee}_w$ are equivalent to generalized cluster 
categories, including $\Cc_Q$ and $\underline{\mod} \Lambda$, where $\Lambda$ is preprojective of Dynkin type. This result is
extended to the case of $c$-sortable words in \cite{AIRT}, with a
similar choice for $\A$. One of the main results in this paper is the following: Each category  $\underline{\Ee}_w$
associated with an element $w$ in a Coxeter group is equivalent to a generalized cluster category $\Cc_{\A}$ for some algebra $\A$
of global dimension at most 2 (Theorem~\ref{applicationbirs}).

 Actually, we prove our main result in a more general setting: We start with a Frobenius category $\Ee$, which we assume to be $\Hom$-finite, stably 2-CY and which has a cluster-tilting object $T$. We assume that the endomorphism algebra $\End_{\underline{\Ee}}(T)$ is Jacobian and has a grading with certain properties. From these data we construct an algebra $\A$ of global dimension at most $2$ and a triangle equivalence $\Cc_{\A}\simeq \underline{\Ee}$ (Theorem~\ref{theoremgrading}) sending the canonical cluster-tilting object $\A$  of $\Cc_{\A}$ to the cluster-tilting object $T$ in $\underline{\Ee}$. The algebra $\A$ is constructed as the degree zero part of $\End_{\underline{\Ee}}(T)$, and we show that $\End_{\underline{\Ee}}(T)$ and $\End_{\Cc_{\A}}(\A)$ are isomorphic algebras (Proposition~\ref{algebraBbar}). This is an important step in the proof of the equivalence. It is however not known in general if 
2-CY categories are equivalent when they have cluster-tilting objects whose endomorphism algebras are isomorphic. The only general result known of this type is that
if the quiver $Q$ of $\End_{\Cc}(T)$ has no oriented cycles, where $T$ is a cluster-tilting object in an algebraic $2$-CY
category $\Cc$, then $\Cc$ is triangle equivalent to the cluster category $\Cc_Q$ \cite{Kel4}. A crucial step in this paper for proving the equivalence $\Cc_{\A}\simeq \underline{\Ee}$ is the construction of a triangle functor from $\Cc_{\A}$ to $\underline{\Ee}$, sending $\A$ to $T$. 
This is done by first constructing a triangle functor from $\Dd^{\rm b}(\A)$ to $\underline{\Ee}$, with strong use of the Frobenius structure of $\Ee$. 

 It is also important to deal with the endomorphism algebra $\End_{\Ee}(T)$ of the cluster-tilting object $T$ in the Frobenius category $\Ee$, rather than only $\underline{\End}_{\Ee}(T)$.  We assume that this algebra is a graded frozen Jacobian algebra (see section 1 for definitions), with potential homogeneous of degree 1. 
Theorem~\ref{theoremgrading} applies in particular to the categories $\Ee_w$ associated to any element in a Coxeter group associated with a finite quiver without oriented cycles (see section 4).

\medskip

The paper is organized as follows. In section 1 we recall some background material on cluster-tilting objects in 2-CY categories, on generalized cluster categories from \cite{Ami3} and on Jacobian algebras from \cite{DWZ}, together with the generalization to frozen Jacobian algebras given in \cite{BIRSm}.

 In section 2 we construct a special triangle (Proposition~\ref{triangle}), which is useful for our construction of a functor from $\Cc_{\A}$ to $\underline{\Ee}$.

Section 3 is devoted to the proof of the triangle equivalence from $\Cc_{\A}$ to $\underline{\Ee}$ (Theorem~\ref{theoremgrading}). We first show that the global dimension of $\A$ is at most 2. Then we construct our triangle functor from $\Cc_{\A}$ to $\underline{\Ee}$ using the special triangle from section 2, together with a universal property from \cite{Kel, Ami3}. Finally we show that our functor is an equivalence by using a criterion from \cite{Kel4}.

In section 4 we apply the main theorem to prove that for any element $w$ in a Coxeter group, the 2-CY triangulated category $\underline{\Ee}_w$ is triangle equivalent to some generalized cluster category, which was our original motivation (Theorem~\ref{applicationbirs}).

In section 5 we give two examples to illustrate our results. The first one is an illustration of Theorem~\ref{applicationbirs}. In the second one we use Theorem~\ref{theoremgrading} to construct a triangle equivalence from a generalized cluster category $\Cc_{\A}$ to a category $\underline{\Ee}_w$ which sends the canonical cluster-tilting object $\A$ to a cluster-tilting object $T$ in $\underline{\Ee}_w$, where $T$ is not associated to a reduced expression of $w$.

\subsection*{Notations}
 By a triangulated category we mean a $k$-linear triangulated category satisfying the Krull-Schmidt property. For all triangulated categories we will denote the shift functor by $[1]$. By a Frobenius category we mean an exact $k$-category with enough projectives and injectives ,where the projectives and the injectives coincide. If $\Ee$ is Frobenius, then the associated stable category $\underline{\Ee}$ is triangulated \cite{Hap88} and is by definition algebraic. For an object $T$ in an additive $k$-category, we denote by $\add(T)$ the additive closure of $T$. For a  $k$-algebra $A$, we denote by $\Mod A$ the category of right $A$-modules and by $\mod A$ the category of finitely presented right $A$-modules. We also denote by $\Dd(A)$ the derived category $\Dd (\Mod A)$ and by $\Dd^{\rm b}A$ the bounded derived category $\Dd^{\rm b}(\mod A)$.  Let $D$ be the usual duality $\Hom_k(?,k)$. The tensor product $-\ten-$, when not specified, is over the ground field $k$. For a quiver $Q$ we denote by $Q_0$ its set of vertices, by $Q_1$ its set of arrows, by $s$ the source map and by $t$ the target map.

\section{Background}

In this section we collect some background material relevant for this paper.

\subsection{Cluster-tilting objects}

Let $\Cc$ be a $k$-category which 
is $\Hom$-finite, that is, has finite dimensional homomorphism spaces over $k$. 
Assume that $\Cc$ is either Frobenius stably 2-CY (that is, its stable category is 2-CY) or triangulated 2-CY. 
Then an object $T$ in $\Cc$ is said to be \emph{cluster-tilting} if 

\hspace{.5cm} (i) $T$ is rigid, \emph{i.e.} $\Ext^1_{\Cc}(T,T)=0$, and

\hspace{.5cm} (ii) $\Ext^1_{\Cc}(T,X)=0$ implies that $X$ is a summand of a finite direct sum of copies of $T$.

Note that when $\Cc$ is Frobenius stably 2-CY, then any indecomposable projective-injective module is a summand of every cluster-tilting object.
The finite dimensional algebras $\End_{\Cc}(T)$, where $\Cc$ is triangulated 2-CY, are called 
\emph{2-CY-tilted algebras}.  

Assume that $T=T_1\oplus\ldots\oplus  T_n$ is a cluster-tilting object in a triangulated 2-CY category $\Cc$,
 where the $T_i$ are indecomposable and pairwise not isomorphic. Then for each $i=1,\ldots,n$ 
there is a unique indecomposable object $T_i^*$ not isomorphic to $T_i$, such that $T^*=(T/T_i)\oplus T_i^*$ 
is a cluster-tilting object \cite{Bua},\cite{Iya}. The new object $T^*$ is called the \emph{mutation} of $T$ at $T_i$.

If $T=T_1\oplus\ldots\oplus  T_n$ is a cluster-tilting object in a Frobenius stably 2-CY category, 
we can only mutate at the $T_i$ which are not projective-injective.

When $T_i^*$ is defined, there are \emph{exchange sequences} if $\Cc$ is Frobenius 
$$\xymatrix{0\ar[r] & T_i^*\ar[r]^f & B\ar[r]^g & T_i\ar[r] & 0 &  \textrm{and} &  
0\ar[r] & T_i\ar[r]^{f'} & B'\ar[r]^{g'} & T_i^*\ar[r] & 0}$$ 
or \emph{exchange triangles} if $\Cc$ is triangulated
$$\xymatrix{ T_i^*\ar[r]^f & B\ar[r]^g & T_i\ar[r] & T_i^*[1] & \textrm{ and} &  
 T_i\ar[r]^{f'} & B'\ar[r]^{g'} & T_i^*\ar[r] & T_i[1]}$$ where $f$, $f'$ are 
minimal left $\add(T/T_i)$-approximations and $g$, $g'$ are minimal right $\add(T/T_i)$-approximations.
These sequences (or triangles) play an important role in the categorification of cluster algebras.
\medskip

There is also a related kind of sequences investigated in \cite{Iya}. 

\begin{prop}\label{2almostsplit}[Iyama-Yoshino]
Let as before $\Cc$ be a $\Hom$-finite Frobenius stably 2-CY category with a cluster-tilting object $T=T_1\oplus\ldots\oplus T_n$. 
For each $i=1,\ldots,n$, 
if $T_i$ is not projective-injective, there are exact sequences 
$$\xymatrix{0\ar[r] & T_i^+\ar[r]^f & E\ar[r]^g & T_i\ar[r] & 0 &  \textrm{and} &  
0\ar[r] & T_i\ar[r]^{f'} & E'\ar[r]^{g'} & T_i^+\ar[r] & 0}$$ for some indecomposable object $T_i^+$ in $\Cc$,
such that $g$ (resp. $g'$) is right almost split in $\add(T)$ (resp. in $\add((T/T_i)\oplus T_i^+)$) and $f'$ (resp. $f$) is left almost split in $\add(T)$ (resp. in $\add((T/T_i)\oplus T_i^+)$). 

\end{prop}
\noindent
The induced sequence 
$\xymatrix{0\ar[r] & T_i\ar[r]^{f'} & E'\ar[r]^{fg'} & E\ar[r]^g & T_i\ar[r] & 0}$ is called the
\emph{2-almost split sequence} associated with $T_i$.
\medskip

There is a corresponding result when $\Cc$ is triangulated. For any indecomposable direct summand $T_i$ of a cluster-tilting object $T$, there are triangles
$$\xymatrix{ T_i^+\ar[r]^f & E\ar[r]^g & T_i\ar[r] & T_i^+[1] & \textrm{ and} &  
 T_i\ar[r]^{f'} & E'\ar[r]^{g'} & T_i^+\ar[r] & T_i[1]})$$
where the maps $f,f'$ are left almost split and $g,g'$ are right almost split. 
For cluster categories, it was shown in \cite{Bua} that these triangles coincide
 with the exchange triangles. 
More generally, they clearly coincide with the exchange triangles if and only if there are no loops in the quiver of $\End_{\Cc}(T)$.

\medskip
Using the existence of 2-almost split sequences in a $\Hom$-finite Frobenius category, we can construct minimal projective and injective resolutions of simple modules over $\End_\Cc(T)$. Note that it already follows from \cite[2.5.3]{Iya07b} that $\End_{\Cc}(T)$ has global dimension at most 3.

\begin{prop}\label{projresolsimple}
Let $\Cc$ be a $\Hom$-finite Frobenius stably 2-CY category with a cluster-tilting object $T=T_1\oplus\ldots\oplus T_n$, where the $T_i$ are indecomposable and pairwise non isomorphic. Let $B:=\End_\Cc(T)$ be the endomorphism algebra of $T$, and let $Q$ be the quiver of $B$, which is then isomorphic to $kQ/I$ for some admissible ideal $I$ in $kQ$. For each $i=1,\ldots,n$, 
such that $T_i$ is not projective-injective, denote by $S_i$  the simple $B$-module $\Hom_\Cc(T,T_i)/Rad(\Hom_\Cc(T,T_i))$.
Then the minimal projective and injective resolutions of the simple $B$-module $S_i$ are of the form:
$$\xymatrix{0\ar[r] & e_iB\ar[r]^(.3){(b)} & \bigoplus_{b\in Q_1,s(b)=i}e_{t(b)}B \ar[r]^{(r_{ab})} & \bigoplus_{a\in Q_1,t(a)=i}e_{s(a)}B\ar[r]^(.7){(a)} & e_iB\ar[r] & S_i\ar[r] & 0}$$
$$\xymatrix@-.2cm{0\ar[r] & S_i\ar[r] & D(Be_i)\ar[r]^(.3){(b)} & \bigoplus_{b\in Q_1,s(b)=i}D(Be_{t(b)}) \ar[r]^{(r'_{ab})} & \bigoplus_{a\in Q_1,t(a)=i}D(Be_{s(a)})\ar[r]^(.7){(a)} & D(Be_i)\ar[r] & 0},$$
for some maps $r_{a,b}$, where the sets $\{ar_{ab}|\  a,b\in Q_1, t(a)=s(b)=i, t(b)=j\}$ and  $\{r'_{ab}b|\  s(b)=t(a)=i, s(a)=j\}$ are bases of $e_iIe_j$. 
\end{prop}

\begin{proof}

Applying the functor $\Hom_{\Cc}(T,-)$ to the 2-almost-split sequence  $$\xymatrix{0\ar[r] & T_i\ar[r]^{f'} & E'\ar[r]^{fg'} & E\ar[r]^g & T_i\ar[r] & 0}$$ we get the following exact sequence of $B$-modules

$$\xymatrix{0\ar[r] & \Hom_\Cc(T,T_i)\ar[r] & \Hom_\Cc(T,E')\ar[r] & \Hom_\Cc(T,E)\ar[r] & \Hom_\Cc(T,T_i)\ar[r] & S_i\ar[r] &  0} $$
which is a minimal projective resolution of the simple $B$-module $S_i$.

Let $Q$ be the quiver of $B$, and $B\simeq kQ/I$. Since $g$ is minimal right almost split in $\add(T)$, we have $$E\simeq \bigoplus_{a\in Q_1|\ t(a)=i}T_{s(a)}\quad  \textrm{and}\quad \Hom_\Cc(T,g)\simeq (a)_{\{a\in Q_1|\ t(a)=i\}}. $$ Since $f'$ is minimal left almost split in $\add(T)$, we have 
$$E'\simeq \bigoplus_{b\in Q_1|\ s(b)=i}T_{t(b)} \quad \textrm{and} \quad \Hom_\Cc(T,f')\simeq (b)_{\{b\in Q_1| \ s(b)=i\}} .$$ 

For $a,b\in Q_1$ with $t(a)=i$ and $s(b)=i$, let $r_{ab}:e_{t(b)}B\rightarrow e_{s(a)}B$ be the map induced by $fg':\bigoplus_{b\in Q_1|\ s(b)=i}T_{t(b)}\rightarrow\bigoplus_{a\in Q_1|\ t(a)=i}T_{s(a)} $. Since the 2-almost split sequence associated to $T_i$ induces a minimal projective resolution of the simple $B$-module $S_i$, the set $\{ar_{ab}|\  a,b\in Q_1, t(a)=s(b)=i, t(b)=j\}$ is a basis of the set of relations $e_iIe_j$.

To get the other sequence of the proposition, we apply the functor $D\Hom_\Cc(-,T)$ to the 2-almost split sequence associated to $T_i$, and we proceed similarly.

\end{proof}

\subsection{Generalized cluster categories}

Let $\Lambda$ be a finite dimensional $k$-algebra of global dimension
at most $2$. We denote by $\Dd^{\rm b}(\Lambda)$ the bounded derived category
of finitely generated (right) $\Lambda$-modules. It has a Serre functor that
we denote by $\mathbb{S}$, which coincides with $\tau[1]$. 

The \emph{generalized cluster category} $\Cc_\Lambda$ is defined in \cite{Ami3} as
the triangulated hull $(\Dd^{\rm b}(\Lambda)/\mathbb{S}[-2])_\Delta$ in the sense of~\cite{Kel} of the orbit category $\Dd^{\rm b}(\Lambda)/\mathbb{S}[-2]$. The composition of the functors 
$$\xymatrix{\pi_\Lambda:\Dd^{\rm b}(\Lambda)\ar@{->>}[r] &
  \Dd^{\rm b}(\Lambda)/\mathbb{S}[-2]\ar@{^(->}[r] & \Cc_\Lambda}$$ is a
triangle functor.

The following definition and theorem give a more explicit description of the generalized cluster category. This is added here for the convenience of the reader but will not be used later.
\begin{dfa}\cite{Kel10}
Denote by $\Theta_2$ a cofibrant resolution of the complex of
$\Lambda$-bimodules $\RHom_\Lambda(D\Lambda,\Lambda)[2]$, that is a complex of projective $\Lambda$-bimodules which is quasi-isomorphic to $\RHom_{\Lambda}(D\Lambda,\Lambda)[2]$. Then the  \emph{derived $3$-preprojective algebra} is
defined as the tensor DG algebra:
\[ \BPi_{3}(\Lambda):= T_\Lambda \Theta_2=\Lambda\oplus \Theta_2\oplus (\Theta_2\ten_\Lambda \Theta_2)\oplus \ldots .\]
We set $\Pi_{3}(\Lambda):=H^0(\BPi_{3}(\Lambda)) $. It is called
the \emph{$3$-preprojective algebra}.
\end{dfa}

\begin{thma}\cite{Kel,Ami3}
Let $\Lambda$ be a finite dimensional algebra of global dimension at most 2. Then there exists a triangle equivalence \[\Cc_\Lambda:=(\Dd^{\rm
b}\Lambda/\mathbb{S}[-2])_\Delta\simeq \per \BPi_{3}(\Lambda)/\Dd^{\rm
b}\BPi_{3}(\Lambda )\]
where $\per \BPi_3(\Lambda)$ is the thick subcategory of $\Dd \BPi_3(\Lambda)$ generated by $\BPi_3(\Lambda)$, and $\Dd^{\rm b}\BPi_3(\Lambda)$ is the thick subcategory of $\Dd\BPi_3(\Lambda)$ formed by the objects having finite dimensional total cohomology.
\end{thma}

There is a useful criterion for constructing triangle functors from a 
generalized cluster category $\Cc_\Lambda$ to some stable category $\underline{\Ee}$. It can be deduced from the 
universal property of $\pi_\Lambda$ given in subsection 4.1 of \cite{Ami3} (see also section 9 of \cite{Kel} or the appendix of \cite{IO2} for more details). This criterion, which is given in the next proposition, is a key step for proving the equivalence of the main theorem of this paper. 

For a Frobenius category $\Ee$ and an algebra $\Lambda$, we here denote by $\Dd^{\rm b}(\Lambda^{\rm op}\ten\Ee)$ the bounded derived category of the exact category $(\mod\Lambda^{\rm op})\ten\Ee$ whose objects are objects in $\Ee$ having a structure of finitely generated left $\Lambda$-module (see for example \cite[sections 2 and 3]{Kel06} for precise definitions of tensor products of $k$-categories and derived categories).

\begin{prop}\label{univprop}
Let $\Cc_\Lambda$ be a generalized cluster category, where $\Lambda$ is a finite dimensional algebra of global dimension at most $2$. Let $\Ee$ be a
Frobenius category. Let $M$ be an object in $\Ee$ and assume that $M$ has a left
$\Lambda$-module structure. Assume that there is a morphism in $\Dd^{\rm b} (\Lambda^{\rm op}\ten \Ee)$
$$\alpha: M\longrightarrow \RHom_\Lambda (D\Lambda,\Lambda)\lten_\Lambda M[2]$$ 
whose cone lies in $\Dd^{\rm b}(\Lambda^{\rm op}\ten \Pp)$, where $\Pp$ is the full 
subcategory of $\Ee$ of projective-injectives. Then there exists
a triangle functor $\Cc_\Lambda\rightarrow \underline{\Ee}$ such that the following diagram commutes
$$\xymatrix{\Dd^{\rm b}(\Lambda)\ar[rr]^{-\lten_\Lambda M}\ar[d]^{\pi_\Lambda} && \Dd^{\rm b}(\Ee)\ar[d] \\ \Cc_\Lambda\ar[rr] && \underline{\Ee}}.$$
\end{prop}

Note that the endofunctor $-\lten_\Lambda \RHom_\Lambda(D\Lambda,\Lambda)[2]\simeq \RHom_\Lambda(D\Lambda,-)[2]$ of $\Dd^{\rm b}(\Lambda)$ is isomorphic to the functor $\mathbb{S}^{-1}[2]$. Hence Proposition~\ref{univprop} requires in particular that for any $X$ in $\Dd^b(\Lambda)$, there is a morphism $X\lten_{\Lambda}M\rightarrow X\lten_{\Lambda}\RHom_{\Lambda}(D\Lambda,\Lambda)\lten_{\Lambda}M[2]\simeq\mathbb{S}^{-1} X[2]\lten_{\Lambda}M$ in $\Dd^b(\Ee)$ whose cone is in $\Dd^b(\Pp)$. In other words,  the images of $X$ and of $\mathbb{S}^{-1}X[2]$  under the composition
$$\xymatrix{ \Dd^{\rm b}(\Lambda)\ar[rr]^{-\lten_\Lambda M} && \Dd^{\rm b}(\Ee)\ar[r] & \Dd^{\rm b}(\Ee)/\Dd^{\rm b}(\Pp)=\underline{\Ee}}$$
are isomorphic. Here the category $\Dd^{\rm b}(\Pp)$ is the thick subcategory of $\Dd^{\rm b}(\Ee)$ generated by $\Pp$. Thus the localization of $\Dd^{\rm b}(\Ee)$ by $\Dd^{\rm b}(\Pp)$ is equivalent to the stable category $\underline{\Ee}$ by \cite{Kel6}, and this localization gives us the right vertical map of this diagram. This implies immediately that the functor $\lten_{\Lambda}M:\Dd^b(\Lambda)\rightarrow \underline{\Ee}$ factors through the orbit category $\Dd^b(\Lambda)/\mathbb{S}[-2]$. However the proof that it factors through the generalized cluster category as a triangle functor is non trivial and uses the universal property of triangulated orbit categories \cite{Kel,IO2}.

\medskip

The case when the 3-preprojective algebra of $\Lambda$ is finite dimensional is especially nice.
 \begin{thma}\cite[Theorem 4.10]{Ami3}
Let $\Lambda$ be a finite dimensional  algebra of global dimension
at most $2$, and assume that the algebra $\Pi_3(\Lambda)$  is finite dimensional.
Then the category $\Cc_\Lambda$ is a $\Hom$-finite 
$2$-Calabi-Yau category, the object $\pi_\Lambda(\Lambda)\in\Cc_\Lambda$ is a cluster-tilting object and we have an isomorphism \[\End_{\Cc_\Lambda}(\pi_\Lambda(\Lambda))\simeq \Pi_3(\Lambda).\]
\end{thma}
\subsection{Jacobian algebras and generalizations}

Quivers with potentials and their associated Jacobian algebras have been investigated in \cite{DWZ}.
Let $Q$ be a finite quiver. For each arrow $a$ in $Q$, the \emph{cyclic derivative} $\partial_a$ with
respect to $a$ is the unique linear map
$\partial_a:kQ\rightarrow kQ$
which takes the class of a path $p$ to the sum $\sum_{p=uav}vu$ taken over all decompositions of the 
path $p$ (where $u$ and $v$ are possibly idempotent elements $e_i$ associated to the vertex $i$).
A \emph{potential} on $Q$ is any linear combination $W$ of cycles 
in $Q$. The associated Jacobian algebra  is by definition the algebra
$$\Jac(Q,W):=kQ/\langle \partial_aW; a\in Q_1\rangle.$$
There is a more general definition given in \cite{DWZ}, 
dealing with the complete path algebras, and hence there is also 
a larger class of Jacobian algebras. However, in this paper we only consider the Jacobian algebras defined above.

Any finite dimensional Jacobian algebra (in the general sense) is 2-CY-tilted (\cite{Ami3,Kel10}). As a partial converse, some classes of 2-CY-tilted algebras associated with elements in Coxeter groups are Jacobian (\cite{BIRSm}). Furthermore, the 2-CY-tilted algebras given by the canonical cluster-tilting object in a generalized cluster category are Jacobian, as stated in the following result \cite[Theorem 6.11 \textit{a)}]{Kel10}. 
 
\begin{thma}[Keller]\label{keller}
Let $A=kQ/I$ be an algebra of global dimension at most $2$, such that $I$
is generated by a finite set of minimal relations $(r_i)$. The relation
$r_i$ starts at the vertex $s(r_i)$ and ends at the vertex
$t(r_i)$. Let $\widetilde{Q}$ be the quiver obtained from $Q$ by adding additional arrows $a_i:t(r_i)\rightarrow s(r_i)$ for each minimal relation $r_i$, and let $W_A$ be the potential $\sum_i a_ir_i$. Then there is an isomorphism of algebras:
$$\End_{\Cc_A}(A)\simeq \Jac(\widetilde{Q},W_A).$$ 
\end{thma} 
\medskip
There is a generalization of quivers with potentials $(Q,W)$ to \emph{frozen quivers with potentials} $(Q,W,F)$ 
in \cite{BIRSm}. Here $F=(F_0,F_1)$ is a pair consisting of a subset $F_0$ of vertices of $Q$ (called \emph{frozen vertices})
and the subset $F_1=\{a\in Q_1, s(a)\in F_0 \textrm{ and }t(a)\in F_0\}$ of arrows (called \emph{frozen arrows}).
The associated \emph{frozen Jacobian algebra} is by definition the algebra
$$\Jac(Q,W,F):=kQ/\langle \partial_aW, a\notin F_1\rangle$$ 

As for ordinary quivers with potential  in \cite{DWZ}, one can define a \emph{reduced} frozen quiver with potential $(Q,W,F)$ by requiring that each term in $W$  has length at least 3 and has at least one arrow in $Q_1\setminus F_1$.

\section{A useful triangle}

In this section we study the Calabi-Yau property of a frozen Jacobian algebra $B$, when $B$ is the endomorphism algebra of a cluster-tilting object in a Frobenius stably 2-Calabi-Yau category (subsection~2.2). Then we assume that $B$ has a special grading and construct some algebras $A$ and $\bar{A}$ related to $B$ (subsection~2.3). Using the grading on $B$ and the Calabi-Yau property we construct a triangle in the category $\Dd^{\rm b}(\A^{\rm op}\ten B)$, which will be crucial for the proof of our main result (subsections~2.4 and~2.5).

\subsection{Basic setup}

Let $\Ee$ be a Frobenius category which is $\Hom$-finite and stably $2$-Calabi-Yau. While we are mainly interested in the stable category $\underline{\Ee}$, it will be important to first consider a cluster-tilting object $T$ in the Frobenius category $\Ee$, and its endomorphism algebra $B:=\End_{\Ee}(T)$. We assume that this algebra $B$ is isomorphic to $\Jac(Q,W,F)$ for some reduced frozen quiver with potential $(Q,W,F)$. Notice that the quiver of $B$ is $Q$ since the potential $W$ is reduced.  We also assume the following:

\medskip
\hspace{.5cm}(H1) The vertices in $F_0$ correspond to the isoclasses of projective-injective indecomposables in the Frobenius category $\Ee$.

\medskip
\hspace{.5cm}(H2) The set  $\{\partial_aW, a\notin F_1\}$, which by definition generates the ideal of relations of the algebra $\Jac(Q,W,F)$, forms a $k$-basis for this ideal. In particular, for any $a\in Q_1\setminus F_1$, we have $\partial_aW\neq 0$.

\medskip
For $i\in Q_0$ we denote by $e_i$ the primitive idempotent of $B$ associated to $i$. Denote by $e_F$ the idempotent $\sum_{i\in F_0} e_i$ and consider the factor algebra $\bar{B}:=B/Be_FB$. By (H1) we have an isomorphism of algebras $\bar{B}\simeq\End_{\underline{\Ee}}(T)$.

Let $\bar{Q}$ be the full subquiver of $Q$ obtained from $Q$ by deleting the vertices in $F_0$. We then have a projection 
$\xymatrix{kQ\ar@{->>}[r] & kQ/kQe_FkQ\simeq k\bar{Q}}.$ We denote by $\bar{W}$ the image of $W$ under this projection.
It is not hard to see that there is an isomorphism $\bar{B}\simeq \Jac(\bar{Q},\bar{W})$. Indeed, since any arrow $a$ in 
$F_1$ satisfies $s(a)\in F_0$ and $t(a)\in F_0$, the partial derivative $\partial_a\bar{W}$ vanishes
for any $a$ in $F_1$.

\subsection{Calabi-Yau property for $B$}
In this subsection we describe the minimal projective and injective resolutions of $\bar{B}$ as $\bar{B}$-$B$-bimodule and deduce a Calabi-Yau property linking these two algebras. We start with giving explicit projective and injective resolutions over $B$ of the simple $\bar{B}$-modules.  
\begin{lema}\label{projresolS}
Let $B\simeq \Jac(Q,W,F)$ be the endomorphism algebra of a cluster-tilting object $T$ in a $\Hom$-finite Frobenius stably 2-CY category $\Ee$. Assume that (H1) and (H2) hold. Then for any $i\in Q_0\setminus F_0$, the sequences
$$\xymatrix{0\ar[r] & e_iB\ar[r]^(.3){(b)} & \bigoplus_{b,s(b)=i}e_{t(b)}B \ar[rr]^{(a^{-1}(\partial_{b}W))} && \bigoplus_{a,t(a)=i}e_{s(a)}B\ar[r]^(.7){(a)} & e_iB\ar[r] & S_i\ar[r] & 0}$$
and
$$\xymatrix@C=.5cm{0\ar[r] & S_i\ar[r] & D(Be_i)\ar[r]^(.3){(b)} & \bigoplus_{b,s(b)=i}D(Be_{t(b)}) \ar[rr]^{((\partial_{a}W)b^{-1})}& & \bigoplus_{a,t(a)=i}D(Be_{s(a)})\ar[r]^(.7){(a)} & D(Be_i)\ar[r] & 0}$$
are minimal projective and injective resolutions of the simple $B$-module $S_i$. When $v$ is a path in $Q$, we write $a^{-1}v=u$ if $v=au$ in $kQ$ and $a^{-1}v=0$ otherwise.
\end{lema}

\begin{proof}
By (H1), if $i$ is not in $F_0$, the corresponding summand $T_i$ of $T$ is not projective-injective. Hence by Proposition~\ref{2almostsplit}, there exists a 2-almost split sequence associated with $T_i$ . Moreover, since the potential $W$ is reduced, then $Q$ is the quiver of $\End_{\Ee}(T)$. Thus by Proposition~\ref{projresolsimple} we obtain a minimal projective $B$-resolution of $S_i$ 
$$\xymatrix{0\ar[r] & e_iB\ar[r]^(.3){(b)} & \bigoplus_{b,s(b)=i}e_{t(b)}B \ar[r]^{(r_{ab})} & \bigoplus_{a,t(a)=i}e_{s(a)}B\ar[r]^(.7){(a)} & e_iB\ar[r] & S_i\ar[r] & 0},$$
where $\{ ar_{ab}|\  t(a)=i, s(b)=i\}$ is a basis for the space of relations with target $i$. By (H2) the set $\{ \partial_bW|\ s(b)=i\}$ is a basis for the same space. Thus a minimal projective resolution of $S_i$ can be written as 
$$\xymatrix{0\ar[r] & e_iB\ar[r]^(.3){(b)} & \bigoplus_{b,s(b)=i}e_{t(b)}B \ar[rr]^{(a^{-1}(\partial_bW))} && \bigoplus_{a,t(a)=i}e_{s(a)}B\ar[r]^(.7){(a)} & e_iB\ar[r] & S_i\ar[r] & 0}$$
which proves the first claim. The proof is similar for the injective resolution.
\end{proof}

We use Lemma~\ref{projresolS} to obtain exact sequences of $\bar{B}$-$B$-bimodules which permit to get the functorial minimal projective and injective resolutions of any $\bar{B}$-module when viewed as a $B$-module. This is inspired by \cite{Boc}.

\begin{prop}\label{projresolB}
There exist  exact sequences in $\mod (\bar{B}^{\rm op}\ten B)$

\[\tag{$a$}\xymatrix@C=.5cm{ 0\ar[r] & \bigoplus_{i\notin F_0}P_{i,i}\ar[r] & 
\bigoplus_{a\notin F_1}P_{s(a),t(a)} \ar[r] &
\bigoplus_{a\notin F_1}P_{t(a),s(a)} \ar[r] &
\bigoplus_{i\notin F_0}P_{i,i}\ar[r] & \bar{B} \ar[r] &  0}
\]

\[\tag{$b$}\xymatrix@C=.6cm{ 0\ar[r] &  \bar{B} \ar[r] & \bigoplus_{i\notin F_0}I_{i,i}\ar[r] & 
\bigoplus_{a\notin F_1}I_{s(a),t(a)} \ar[r] &
\bigoplus_{a\notin F_1}I_{t(a),s(a)} \ar[r] &
\bigoplus_{i\notin F_0}I_{i,i}\ar[r] &   0}
\]
where $P_{i,j}:=\bar{B} e_i
\ten e_j B$ and $I_{i,j}:=\Hom_k(Be_j,\bar{B} e_i)$ for $(i,j)\in Q_0\times Q_0$. 
\end{prop}

\begin{rema}\label{remark projective injective}
For $M\in \mod \bar{B}$, if we apply the functor $M\ten_{\bar{B}}-$ to the sequence $(a)$ we get an exact sequence in $\mod B$ since  $P_{i,j}$ is projective as left $\bar{B}$-module. This sequence  is the minimal projective resolution of $M$ when viewed as a $B$-module since $P_{i,j}$ is projective as right $B$-module.

Similarly, if we apply the functor $M\ten_{\bar{B}}-$ to the sequence $(b)$ we get an exact sequence in $\mod B$ since  $I_{i,j}$ is projective as left $\bar{B}$-module. This sequence  is the minimal injective resolution of $M$ when viewed as a $B$-module since $I_{i,j}$ is injective as right $B$-module.
\end{rema}

\begin{proof}
For $(i,j)\in Q_0\times Q_0$ denote by $\Pi_{i,j}$ the projective $B$-bimodule $Be_i\ten e_j
B$. Then consider the following sequence  
\[ \tag{$c$}\xymatrix{\bigoplus_{i\notin F_0}\Pi_{i,i}\ar[r]^-{d_{2}} & 
\bigoplus_{a\notin F_1}\Pi_{s(a),t(a)} \ar[r]^-{d_{1}} &
\bigoplus_{a\in Q_1}\Pi_{t(a),s(a)} \ar[r]^-{d_{0}} &
\bigoplus_{i\in Q_0}\Pi_{i,i}.}
\]
where the maps $d_0$, $d_1$ and $d_2$ are defined as follows: 
\[ \begin{array}{rcl} d_2(e_i\ten e_i) & = & \sum_{a,t(a)=i}a\ten e_i -\sum_{b, s(b)=i}e_i\ten
b; \\
d_1(e_{s(a)}\ten e_{t(a)}) & = & \sum_{b\in Q_1}\partial_{a,b}W  \textrm{ where }\partial _{a,b}(apbq)=p\ten
q\in Be_{t(b)}\ten e_{s(b)}B\\  & & \hspace{2.3cm}\textrm{ for  a cycle } apbq \textrm{ in } Q  ;\\
d_0 (e_{t(a)}\ten e _{s(a)}) & = & a\ten e_{s(a)} -e_{t(a)}\ten a.  \end{array}\]
It is easy to check that this is a complex of $B$-bimodules, and that
$\Coker d_0 =B$. (The map $\bigoplus_{i\in Q_0} \Pi_{i,i}\rightarrow B$ is the
multiplication map.)

Applying the functor $\bar{B}\ten_B -$, we get the complex of
$\bar{B}$-$B$-bimodules 
$$\xymatrix{ \bigoplus_{i\notin F_0}P_{i,i}\ar[r]^-{d_{2}} & 
\bigoplus_{a\notin F_1}P_{s(a),t(a)} \ar[r]^-{d_{1}} &
\bigoplus_{a\notin F_1}P_{t(a),s(a)} \ar[r]^-{d_{0}} &
\bigoplus_{i\notin F_0}P_{i,i}.}
$$
Indeed, if $i\in F_0$, then $\bar{B} e_i\ten e_i B=0$, and if $a\in F_1$, then $\bar{B} e_{t(a)}\ten e_{s(a)} B=0$. 

Applying the functor $S_i\ten_{\bar{B}}-$ with $i\notin F_0$, we get the complex 
$$ \xymatrix{0\ar[r] & e_iB \ar[r]^(.3){(b)} & \bigoplus_{b, s(b)=i}e_{t(b)}B \ar[rr]^{((\partial_aW)b^{-1})}&&
  \bigoplus _{a, t(a)=i} e_{s(a)}B \ar[r]^(.7){(a)} & e_i B}$$
which is exactly the minimal projective resolution of $S_i$ as $B$-module described in Lemma~\ref{projresolS}. Thus the
following sequence in $\mod (\bar{B}^{\rm op}\ten B)$ is exact:  
$$\xymatrix{0\ar[r] & \bigoplus_{i\notin F_0}P_{i,i}\ar[r] & 
\bigoplus_{a\notin F_1}P_{s(a),t(a)} \ar[r] &
\bigoplus_{a\notin F_1}P_{t(a),s(a)} \ar[r] &
\bigoplus_{i\notin F_0}P_{i,i}\ar[r] & \bar{B} \ar[r] &  0,}
$$
so that ($a$) follows.

\smallskip
We now prove the existence of the sequence $(b)$, where the proof is dual. 
Let us define the sequence of $B$-bimodules 
$$\xymatrix{ \bigoplus_{i\in Q_0}\Upsilon_{i,i}\ar[r]^-{d^{0}} & 
\bigoplus_{a\in Q_1}\Upsilon_{s(a),t(a)} \ar[r]^-{d^{1}} &
\bigoplus_{a\notin F_1}\Upsilon_{t(a),s(a)} \ar[r]^-{d^2} &
\bigoplus_{i\notin F_0}\Upsilon_{i,i}}
$$ where $\Upsilon_{i,j}:=\Hom_k(Be_j,
Be_i)$ for $(i,j)\in Q_0\times Q_0$. The maps $d^0, d^1, d^2$ are the following
$$\begin{array}{rcl}
d^0(\phi_i) & = & (\sum_{s(a)=i}a\ten e_i-\sum_{t(b)=i}e_i\ten b).(\phi_i) \\
d^1(\phi_a) &= & \sum_{b\in Q_1} (\partial_{a,b}W).(\phi_a)\\
d^2(\phi_a) & = & (a\ten e_i-e_i\ten a).(\phi_a) \end{array} $$
where $(a\ten b)(\phi)(-)=\phi(-a)b$. For instance we have
$$d^0(\phi_i)(-)=\sum_{s(a)=i}\phi_i(-a)-\sum_{t(b)=i}\phi_i(-)b.$$
The kernel of $d^0$ is $B$, the bimodule map $B\rightarrow \bigoplus
\Upsilon_{i,i}$ maps $1_B$ to $(\bold{1}_{Be_i})_i$. 
Using the fact that $S_l\ten_{\bar{B}} \Hom_k(Be_i, \bar{B} e_j)\simeq
\delta_{j,l} e_i DB$ we get the exact sequence $(b)$. 
\end{proof}

We have the following direct consequence (see also \cite[5.4,Thm (b)]{Kel5}) which we include even though it will not be used later in this paper.
\begin{cora}\label{corCY}
There is an isomorphism 
$\RHom_B(DB,\bar{B})[3]\simeq \bar{B} \textrm{ in } \Dd^{\rm b} (\bar{B}^{\rm op}\ten B).$
\end{cora}
\begin{proof}
Since $I_{i,j}=\Hom_k(Be_j,e_i\bar{B})$ is injective as a right $B$-module, then using the exact sequence~$(b)$ of Proposition~\ref{projresolB},
we obtain that $\RHom_B(DB,\bar{B})$ is isomorphic in $\Dd^{\rm b} (\bar{B}^{\rm op}\ten B)$ to the complex
$$ \xymatrix{\bigoplus_{i\notin F_0}R_{i,i}\ar[r]
  &\bigoplus_{a\notin F_1}R_{s(a),t(a)} \ar[r] & 
\bigoplus_{a\notin F_1}R_{t(a),s(a)} \ar[r] &
\bigoplus_{i\notin F_0}R_{i,i}}
$$
where $R_{i,j}$ is the bimodule $\Hom_B(DB,I_{i,j})$.  
Moreover we have the following isomorphisms 
$$\begin{array}{rcl}\Hom_B(DB,I_{i,j})& = & \Hom_B(DB, \Hom_k(B e_j,
  \bar{B} e_i))\\
 &\simeq & \Hom_k (DB e_j,\bar{B} e_i)\\
 & \simeq & \bar{B} e_i \ten e_j B\\ & = & P_{i,j}
\end{array}$$
in $\mod (\bar{B}^{\rm op}\ten B)$. Then the claim follows from the exact sequence $(a)$ of Proposition~\ref{projresolB}.

\end{proof}

\begin{rema}
If there exists an algebra $B=\Jac(Q,W,F)$ with $F_0=\emptyset$ such that the sequences in Lemma~\ref{projresolS} are exact, then we get an isomorphism in $\Dd(B^{\rm op}\ten B)$ between $\RHom(DB,B)[3]$ and $B$ (since in this case $\bar{B}=B$). This means by definition that $B$ is bimodule 3-CY. This is exactly what Bocklandt proved in  \cite[Theorem 4.3]{Boc}, namely that a Jacobian algebra $B$ is bimodule 3-CY if the sequence $(c)$ is exact. In our setup,  $F_0$ is never empty since it corresponds to the projective-injectives. Also the algebra $B$ is not Jacobian in general, and the Jacobian algebra $\bar{B}$ is never of global dimension 3 (indeed it is either hereditary or of infinite global dimension \cite[section 2 Corollary]{Kel5}). However we have an isomorphism between $\B$ and $\RHom_B(DB,\B)[3]$ in $\Dd^{\rm b}(\B^{\rm op}\ten B)$. 
\end{rema} 

\begin{rema}
In the next subsections, we will assume that moreover the algebra $B$ has a special grading that induces an isomorphism $\RHom_{B}(DB,\B)[3](-1)\simeq \B$ in $\Dd^{\rm b}(\gr (\B^{\rm op}\ten B))$, where $\gr( \B^{\rm op}\ten B)$ is the category of finite dimensional graded $\B^{\rm op}\ten B$-modules and where $(1)$ is the degree-shift in the category of graded $B$-modules. However, since Proposition~\ref{projresolB} and Corollary~\ref{corCY} hold without any grading hypothesis, we have separated subsection 2.2 from subsections 2.3, 2.4 and 2.5.
 \end{rema}

\subsection{Construction of the algebras $A$ and $\A$}
In order to identify appropriate subalgebras of $B$ which should give rise to the generalized cluster categories we are looking for, it will be convenient to introduce some special gradings on the quiver. Assume as before that $B=\End_{\Ee}(T)$ is isomorphic to some frozen Jacobian algebra $\Jac(Q,W,F)$.
Then we assume that there exists a degree map $\varphi: Q_1\rightarrow \{0,1\}$ with the following property:
 
\medskip
\hspace{.5cm} (H3) The potential $W$ is homogeneous of degree 1.

\medskip
Since the potential is homogeneous, any relation $\partial_aW$
is homogeneous. Hence $\varphi$ induces a grading on $B$. Define the algebras $A$ and $\bar{A}$ by $A:=B_0\quad \textrm{and} \quad \bar{A}:=A/Ae_F A.$ We have surjective algebra maps: $B\rightarrow A\rightarrow \A$. We want to show the following, which is the main result of the section. 
\begin{prop}\label{triangle}
Let $\Ee$ be a Frobenius category which is $\Hom$-finite and stably $2$-Calabi-Yau. We assume that there exists a cluster-tilting object $T$ in $\Ee$ such that its endomorphism algebra is isomorphic to $\Jac(Q,W,F)$ for some reduced frozen quiver with potential $(Q,W,F)$. With the assumptions (H1), (H2) and (H3), there exists a triangle
\[\tag{$*$} \xymatrix{\RHom_A(DA,\bar{A})\lten_A B[2]\ar@<-.5ex>[r] & \A\lten_A
  B\ar@<-.5ex>[r] & \A\ar@<-.5ex>[r] & \RHom_A(DA,\bar{A})\lten_A B[3]} \] 
in $\Dd^{\rm b}(\A^{\rm op}\otimes B)$, where $A=B_0$ and $\A=A/Ae_FA$.
\end{prop}

The proof is given in the next subsections.  In subsection 2.4 we construct a triangle $X\rightarrow Y\rightarrow \A\rightarrow X[1]$ in $\Dd^{\rm b}(\A^{\rm op}\ten B)$. In subsection 2.5, we show that this triangle is the triangle~$(*)$ of Proposition~\ref{triangle}.

\subsection{A triangle $X\rightarrow Y\rightarrow \A\rightarrow X[1]$}

The following  is proved the same way as Proposition~\ref{projresolB}. We describe exact sequences of $\A$-$B$-bimodules, which give the functorial minimal projective and injective resolutions of the $\A$-modules when viewed as $B$-modules (cf Remark~\ref{remark projective injective}).

\begin{prop}\label{projresol}
There exist  exact sequences in $\mod (\A^{\rm op}\ten B)$

\[ \tag{$a'$}\xymatrix@C=.5cm{0\ar[r] & \bigoplus_{i\notin F_0}P_{i,i}\ar[r]^(.4){d_2} & 
\bigoplus_{a\notin F_1}P_{s(a),t(a)} \ar[r]^{d_1} &
\bigoplus_{a\notin F_1}P_{t(a),s(a)} \ar[r]^(.6){d_0} &
\bigoplus_{i\notin F_0}P_{i,i}\ar[r] & \A \ar[r] &  0}
\]

\[\tag{$b'$} \xymatrix@C=.6cm{0\ar[r] &  \A \ar[r] & \bigoplus_{i\notin F_0}I_{i,i}\ar[r] & 
\bigoplus_{a\notin F_1}I_{s(a),t(a)} \ar[r] &
\bigoplus_{a\notin F_1}I_{t(a),s(a)} \ar[r] &
\bigoplus_{i\notin F_0}I_{i,i}\ar[r] &   0}
\]
where $P_{i,j}:=\A e_i
\ten e_j B$, $I_{i,j}:=\Hom_k(Be_j,\A e_i)$ and the $d_i$ are defined as in Proposition~\ref{projresolB}. 
\end{prop}

Let $a\in Q_1\setminus F_1$ be an arrow with $\varphi(a)=1$. By definition
$d_1(e_{s(a)}\ten e_{t(a)})=\sum_{b\in Q_1}\partial_{a,b}W$.
If $apbq$ is a cycle in $W$, then $\partial_{a,b}(apbq)=p\ten q$ lies
in $\A e_{t(b)}\ten e_{s(b)}B$. Since the degree of $a$ is 1, and since
we have assumed that the potential $W$ is homogeneous of degree 1, then $b$ is of degree
0. Hence the restriction of $d_1$ 
$$\xymatrix{\bigoplus_{a\notin F_1, \varphi(a)=1}P_{s(a),t(a)}\ar[r] & \bigoplus_{a\notin F_1, \varphi(a)=1}P_{t(a),s(a)}}$$
is zero.
Therefore the complex 
$$\xymatrix{\bigoplus_{i\notin F_0}P_{i,i}\ar[r] & 
\bigoplus_{a\notin F_1}P_{s(a),t(a)} \ar[r] &
\bigoplus_{a\notin F_1}P_{t(a),s(a)} \ar[r] &
\bigoplus_{i\notin F_0}P_{i,i}}$$ is isomorphic to the mapping cone of a complex morphism 
 
$$\xymatrix{X:=\bigoplus_{i\notin F_0}P_{i,i}\ar[r]\ar[d] & 
\bigoplus_{a\notin F_1, \varphi(a)=0}P_{s(a),t(a)} \ar[r]\ar[d] &
\bigoplus_{a\notin F_1, \varphi(a)=1}P_{t(a),s(a)}\ar[d] \\
Y:=\bigoplus_{a\notin F_1, \varphi(a)=1}P_{s(a),t(a)} \ar[r]
&\bigoplus_{a\notin F_1, \varphi(a)=0}P_{t(a),s(a)}\ar[r] &
\bigoplus_{i\notin F_0}P_{i,i}.}$$
It is not hard to check that all horizontal maps are homogeneous of
degree 0, and all vertical maps are homogeneous of degree 1. Denote by $f$ the map $X\rightarrow Y$. By Proposition~\ref{projresol}~($a'$), we obtain the triangle 
$\xymatrix{X\ar[r]^f & Y\ar[r] & \A\ar[r] & X[1]}$ in $\Dd(\A^{\rm op}\otimes B).$

\bigskip
Dually, using the exact sequence ($b'$) in Proposition~\ref{projresol},
it is possible to view $\A[1]$ as the mapping cone of a morphism $g:X'\rightarrow Y'$, where horizontal maps are of degree 0, and vertical maps are of degree
$-1$.
$$\xymatrix{X':=\bigoplus_{i\notin F_0}I_{i,i}\ar[r]\ar[d] & 
\bigoplus_{a\notin F_1, \varphi(a)=0}I_{s(a),t(a)} \ar[r]\ar[d] &
\bigoplus_{a\notin F_1, \varphi(a)=1}I_{t(a),s(a)}\ar[d] \\
Y':=\bigoplus_{a\notin F_1, \varphi(a)=1}I_{s(a),t(a)} \ar[r]
&\bigoplus_{a\notin F_1, \varphi(a)=0}I_{t(a),s(a)}\ar[r] &
\bigoplus_{i\notin F_0}I_{i,i}}.$$
 Thus we get a triangle  
$\xymatrix{\A\ar[r] & X'\ar[r]^g & Y'\ar[r] & \A[1]}$ in $\Dd(\A^{\rm op}\otimes B).$

\subsection{Interpretation of $X$ and $Y$}
The aim of this subsection is to construct isomorphisms
$$ Y\simeq \A\lten_A B \quad\textrm{ and }\quad X\simeq \RHom_B(DA,\A)\lten_A
B[2]\quad \textrm{ in } \Dd^{\rm b}(\A^{\rm op}\ten B),$$ in order to prove Proposition~\ref{triangle}. 

We will first show the following
\begin{lema}\label{X_0 and Y_0}
In the setup of subsection~2.4, we have isomorphisms \[ Y_0\simeq \A\quad \textrm{ and} \quad X'_0\simeq \A \quad \textrm{ in } \Dd^{\rm b}(\A^{\rm op}\ten A)\]
\end{lema}

The proof of this result uses the next lemma.  
\begin{lema}\label{lemma graded mapping cone}
Let $B$ be a $\mathbb{Z}$-graded algebra. Let $(X,d_X)$ and $(Y,d_Y)$ be complexes of graded $B$-modules such that the differentials $d_X$ and $d_Y$ are homogeneous of degree $0$. We denote by $(X_p,d_{Xp})$ (resp. $(Y_p,d_{Yp})$) the part of degree $p$ of the complex $(X,d_X)$ (resp. $(Y,d_Y)$). They are complexes of $B_0$-modules. Let $f:X\rightarrow Y$ be a morphism homogeneous of degree $d$, that is the part $f_p$ of degree $p$ of $f$ is a morphism of $B_0$-complexes $X_{p}\rightarrow Y_{p+d}$. Denote by $Z={\sf Cone} (f)$ the mapping cone of $f:X\rightarrow Y$. Then for any integers $p,q$ we have an isomorphism of $B_0$-modules:
$$H^q(Z)_p\simeq H^q ({\sf Cone} ( f_{p-d}:X_{p-d}\rightarrow Y_p)).$$  
\end{lema}

\begin{proof}
The mapping cone of a morphism $f$ homogeneous of degree $d$ in the category of complexes of graded modules with differential homogeneous of degree 0 is still a complex of graded modules with differential homogeneous of degree $0$, and we have \[ {\sf Cone}(f)\simeq \bigoplus_{p\in\mathbb{Z}}{\sf Cone}(f_{p-d}:X_{p-d}\rightarrow Y_p).\]   Then one can check the isomorphisms
\[ H^q({\sf Cone} (f))_p\simeq H^q(({\sf Cone}(f))_p)\simeq H^q({\sf Cone}(f_{p-d})).\]
\end{proof}

\begin{proof}[Proof of Lemma~\ref{X_0 and Y_0}]
Applying Lemma~\ref{lemma graded mapping cone} to the morphism $f:X\rightarrow Y$ defined in subsection~2.4, we get an isomorphism of $(\A^{\rm op}\ten A)$-modules (remember that $A=B_0$):
$$ H^q({\sf Cone} (f:X\rightarrow Y))_0\simeq H^q({\sf Cone}(f_{-1}:X_{-1}\rightarrow Y_0)).$$
By the sequence $(a')$ of Proposition~\ref{projresol} the left term is zero unless $q$ is $0$, and when $q$ is 0, it is isomorphic to $\A$. Since $X$ is non zero only in positive degrees, the right hand side is just $H^q(Y_0)$. Thus we get an isomorphism  $Y_0\simeq \A$ in $\Dd^{\rm b}(\A^{\rm op}\ten A)$.
 
\noindent
Using the triangle $\xymatrix{\A\ar[r] &X'\ar[r]^g & Y'\ar[r] & \A[1] }$ similarly, we get an isomorphism $X'_0\simeq \A$ in $\Dd(\A^{\rm op}\ten A)$.
\end{proof}

The complex $Y_0$ is a complex of  projective $(\A^{\rm op}\ten A)$-modules. Thus for any $\A$-module $M$, $M\lten_{\A} Y_0$ is a complex of projective $A$-modules. By   Lemma~\ref{X_0 and Y_0}, it is quasi-isomorphic to $M$ viewed as an $A$-module. Hence we have the following.
\begin{cora}\label{projdim} Any $\A$-module has projective dimension at most 2 when viewed as an $A$-module.  
\end{cora}

We can now prove our desired isomorphisms.

\begin{lema}\label{interpretation X and Y}
For complexes $X$ and $Y$ defined as in subsection~2.3, there are isomorphisms 
$$ Y\simeq \A\lten_A B \quad\textrm{ and }\quad X\simeq \RHom_B(DA,\A)\lten_A
B[2]\quad \textrm{ in } \Dd^{\rm b}(\A^{\rm op}\ten
 B).$$
\end{lema}

 For the proof of this lemma, we need some basic results:

\begin{lema}\label{lemma1}
Let $B$ be a $\mathbb{Z}$-graded algebra and $A:=B_0$.
Let $P=(P^j,d)$ be a complex of graded projective $B$-modules such that the
differential $d$ of $P$ is homogeneous of degree $0$. Let $P_0$ be the degree 0 part of $P$. Then $P_0$ is a complex over $A$, and we have an
isomorphism of complexes
$$P\simeq P_0\ten_A B. $$
\end{lema}

\begin{proof}
Let $i:A\rightarrow B$ and $p:B\rightarrow A$ be the canonical algebra
maps. We get induced functors:
$$i^*=-\ten_AB: \mod A \rightarrow \mod B, \textrm{ and } p^*=-\ten_B
A: \mod B\rightarrow \mod A$$
Since $B_0\ten_A B\simeq B$, then we have $P^j_0\ten_A B=i^*\circ p^*(P^j)\simeq
P^j$ since $P^j$ is a graded projective $B$-module. Hence, since $P$ is a complex of graded projective modules, we get 
$$P_0\ten_A B =  (P^j_0\ten_A B, i^*\circ p^*
  (d)).$$
Since $i\circ p(b) =b$ if and only if $b$ is of degree $0$, we get $i^*\circ p^*(d)=d$.
\end{proof}

\begin{lema}\label{lemma2}
The functors $\left(\Hom_B(DB,-)\right)_0$ and $\Hom_A(DA,(-)_0)$ are
isomorphic as functors from injective $B$-modules to projective
$A$-modules.
\end{lema}
\begin{proof}
We have $\Hom_{B}(DB,DB)_0\simeq B_0=A$ and $\Hom_{A}(DA,(DB)_0)\simeq \Hom_A(DA,DA)\simeq A$. The rest is easy to check.
It is enough to check it on $DB$, and this is clearly true.
\end{proof} 

\begin{proof}[Proof of Lemma~\ref{interpretation X and Y}]

Since $Y$ is a complex of projective modules, we can apply Lemmas~\ref{lemma1} and~\ref{X_0 and Y_0} to get an isomorphism 
\[ Y\simeq Y_0\ten_A B\simeq \A\lten_A B\quad \textrm{ in } \Dd(\A^{\rm op}\ten
 B).\]

 Using the fact that $\Hom_B(DB,I_{i,j})\simeq P_{i,j}$, we get an isomorphism 
\[ \RHom_B(DB, X') [2]\simeq X \quad \textrm{ in }\Dd(\A^{\rm op}\ten B). \]
Hence we obtain the following isomorphisms in $\Dd(\A^{\rm op}\ten
 B)$
$$\begin{array}{rcll} X & \simeq & \Hom_B(DB,X')[2] & \\
 & \simeq & \left( \Hom_B(DB,X')\right)_0\ten_A B[2] & \textrm{by Lemma
~\ref{lemma1}} \\
& \simeq & \Hom_A(DA,X'_0)\ten_A B[2] & \textrm{by Lemma
 ~\ref{lemma2}}\\
& \simeq & \RHom_A(DA,\A)\lten_A B[2] & \textrm{by Lemma~\ref{X_0 and Y_0}.}
\end{array}$$
\end{proof}

Proposition~\ref{triangle} is a direct consequence of the construction in subsection~2.4 of the triangle $\xymatrix{X\ar[r]^f & Y\ar[r] & \A\ar[r] & X[1]}$ in $\Dd(\A^{\rm op}\otimes B)$ and of Lemma~\ref{interpretation X and Y}.

\section{Main Theorem}

As in the previous section, $\Ee$ is a Frobenius category which is $\Hom$-finite and stably $2$-Calabi-Yau. We assume that there exists a cluster-tilting object $T$ in $\Ee$ whose endomorphism algebra is isomorphic to $\Jac(Q,W,F)$ for some reduced frozen quiver with potential $(Q,W)$ such that (H1), (H2) and (H3) are satisfied. 
Under an additional assumption (H4), we show in this section that the stable category $\underline{\Ee}$ is triangle equivalent to a generalized cluster category.

\subsection{Statement of the main result} In addition to the above assumptions, we assume that we have the following, where $\varphi$ is the degree map required for (H3):

\medskip
\hspace{.5cm} (H4) If $a:i\rightarrow j$ is in $Q_1$ with $i\notin F_0$ and $j\in
  F_0$ then $\varphi(a)=1$.

\medskip
As before we define the algebras $A$ and $\A$ as $A:=B_0\subset B$ and $\A:=A/Ae_FA,$ where $e_F$ is the idempotent $\bigoplus_{i\in F_0}e_i$.
The aim of this section is to prove the following. 
\begin{thma}\label{theoremgrading}
Let $\Ee$ be a Frobenius category which is $\Hom$-finite and stably 2-CY. Assume that there exists a cluster-tilting object $T$ in $\Ee$ such that $B:=\End_\Ee(T)$ is isomorphic to some graded frozen Jacobian algebra $\Jac(Q,W,F)$ which satisfies the conditions (H1)-(H4). Let $A:=B_0$ be the subalgebra of degree 0 of $B$ and $\A:=A/Ae_FA$.  Then we have the following 
\begin{itemize}
\item[($a$)] The algebra $\A$ is of global dimension at most $2$; 
\item[($b$)] There exists a triangle equivalence $\Cc_{\A} \simeq \underline{\Ee},$
where $\Cc_{\A}$ is the generalized cluster category associated to the algebra $\A$.
\end{itemize}
\end{thma}

The proof of the theorem is given in the next three subsections. The fact that $gl.dim \A\leq 2$ is proved in subsection 3.2, the existence of a triangle functor $G:\Cc_{\A}\rightarrow \underline{\Ee}$ is proved in subsection~3.3. Finally it is proved that $G$ is an equivalence of triangulated categories in subsection 3.4.

\subsection{Global dimension of $\A$}
We start with describing the restriction functor $\xymatrix@-.3cm{R:\Dd^{\rm b}(\A)\ar[r]& \Dd^{\rm b}(A)}$ induced by the projection $\xymatrix{A\ar[r] &\A=A/Ae_F A}$. 
\begin{lema}\label{restriction}
Let $\Ee$, $\Jac(Q,W,F)$, $A$ and $\A$ be as in Theorem~\ref{theoremgrading}. Assume that $\Jac(Q,W,F)$ satisfies conditions (H1)-(H4). 
 Then for any $i\notin F_0$, we have the following:
\begin{itemize}
\item[(a)] $R(e_iD\A)\simeq e_iDA$, 
\item[(b)] $R(e_i\A)\simeq \xymatrix{(P\ar[r] & Q\ar[r] & e_i A)}$ where $P$ and $Q$ are in $\add(e_F A)$.
\item[(c)] The functor $R$ is fully faithful.
\end{itemize}
\end{lema}

\begin{proof}
Part $(a)$ follows directly from (H4). 
By Corollary~\ref{projdim} any $\A$-module has projective dimension at most 2 when viewed as an $A$-module. Then $(b)$ follows from (H4). Part $(c)$ follows directly.
\end{proof}

\begin{prop}Let $\Ee$, $\Jac(Q,W,F)$ and $\A$ be as in Theorem~\ref{theoremgrading}, and assume that $\Jac(Q,W,F)$ satisfies  (H1)-(H4).
 Then the global dimension of $\A$ is at most 2.
\end{prop}

\begin{proof}

The complex $X'$, defined in subsection~2.4 by 
\[ \xymatrix{X':= (\bigoplus_{i\notin F_0}I_{i,i}\ar[r]
& \bigoplus_{a\notin F_1, \varphi(a)=0}I_{s(a), t(a)}\ar[r] & \bigoplus_{a\notin F_1, \varphi(a)=1}I_{t(a),s(a)}}), \]
is a complex of graded $(\A^{\rm op}\ten B)$-modules with differential homogeneous of degree 0. The part $X'_0$ of degree 0 is  the complex
\[ \xymatrix{X'_0= (\bigoplus_{i\notin F_0}J_{i,i}\ar[r]
& \bigoplus_{a\notin F_1, \varphi(a)=0}J_{s(a), t(a)}\ar[r] & \bigoplus_{a\notin F_1, \varphi(a)=1}J_{t(a),s(a)}}) \]
where $J_{(i,j)}:=(I_{(i,j)})_0=\Hom_{k}(Be_j,\A e_i)_0\simeq \Hom_k(Ae_j, \A e_i)$.
By Lemma~\ref{X_0 and Y_0}, the complex $X'_0$ is quasi-isomorphic to $\A$. Hence there exists an exact sequence in $\mod (\A^{\rm op}\ten
A)$
$$\xymatrix{0\ar[r] & \A\ar[r] & \bigoplus_{i\notin F_0}J_{i,i}\ar[r]
& \bigoplus_{a\notin F_1, \varphi(a)=0}J_{s(a), t(a)}\ar[r] & \bigoplus_{a\notin F_1, \varphi(a)=1}J_{t(a),s(a)}\ar[r] & 0}$$

Since $J_{i,j}$ is projective as left $\A$-module and injective as right $A$-module, any right $\A$-module, when viewed as an $A$-module, has injective dimension~2 (cf Remark~\ref{remark projective injective}).  By Lemma~\ref{restriction} $(a)$ the injective  $\A$-modules are injective when viewed as $A$-modules. Thus the global dimension of $\A$ is at most 2.
 
\end{proof}

\subsection{Construction of a triangle functor}

Recall that $A=B_0$ is the subalgebra of degree~$0$ of $B=\End_\Ee(T)$. Thus $T$ has a left $A$-module structure. 
Hence we have the following diagram:
$$\xymatrix{\Dd^{\rm b}(\A)\ar[r]_{R}\ar[d]^{\pi_{\A}}\ar@/^.5cm/[rrr]^{-\lten_{\A}(\A\lten_A T)} & \Dd^{\rm b}(A)\ar[rr]_{-\lten_A
    T}& &
  \Dd^{\rm b}(\Ee) \ar[d]  \\ \Cc_{\A}\ar@{..>}[rrr]^{G}&&&\underline{\Ee}.}$$

In this subsection, we prove the existence of a triangle functor $G:\Cc_{\A}\rightarrow \underline{\Ee}$  making the above diagram commute.  This follows from  Proposition~\ref{univprop} (for $M=\A\lten_A T$) together with the following.

\begin{prop}\label{construction alpha}
In the setup of Theorem~\ref{theoremgrading}, there exists a morphism 
\[ \xymatrix{\alpha:\A\lten_A T\ar@<-.5ex>[r] & \RHom_{\A}(D\A,\A)\lten_{\A}(\A\lten_A T)[2]} \quad \textrm{ in } \Dd^{\rm b}(\A^{\rm op}\ten \Ee)\] whose cone is in $\Dd^{\rm b}(\A^{\rm op}\ten \Pp)$, where $\Pp$ is the subcategory of $\Ee$ consisting of the projective-injectives.
\end{prop}

We divide the proof of this result into three lemmas.

\begin{lema}\label{lem1}
There exists an isomorphism \[\xymatrix{\RHom_{A}(DA,\A)\lten_A T[2]\ar@<-.5ex>[r]^-{\sim} & \A\lten_A T} \quad \textrm{ in }\Dd^{\rm b}(\A^{\rm op}\ten \Ee).\]
\end{lema}
\begin{proof}
Applying $-\lten_BT$ to the triangle $(*)$ in Proposition~\ref{triangle}, we
get the following triangle in $\Dd(\A^{\rm op}\ten\Ee)$:
$$\xymatrix{\RHom_A(DA,\bar{A})\lten_A T[2]\ar@<-.5ex>[r] & \A\lten_A
  T\ar@<-.5ex>[r] & \A\lten_BT\ar@<-.5ex>[r] & \RHom_A(DA,\bar{A})\lten_A T[3]}.$$

Using the projective resolution from Lemma~\ref{projresolS} of the simple $B$-module $S_i$ for $i\in Q_0\setminus F_0$ , we conclude that the object $S_i\lten_B T$ is quasi-isomorphic to the complex
$$ \xymatrix@-.8pc{0\ar[r] & e_iB\ten_B T \ar[r] & \bigoplus_{a,
    s(a)=i}e_{t(a)}B\ten_B T \ar[r]&
  \bigoplus _{a, t(a)=i} e_{s(a)}B\ten_B T \ar[r] & e_i B\ten_B T,}$$
which is the 2-almost split sequence associated with $T_i$
$$ \xymatrix{0\ar[r] & T_i \ar[r] & \bigoplus_{a, s(a)=i}T_{t(a)} \ar[r]&
  \bigoplus _{a, t(a)=i} T_{s(a)} \ar[r] & T_i\ar[r] & 0 }.$$
Hence $S_i\lten_B T$ is zero in $\Dd^{\rm b}(\Ee)$. Therefore  for each $M\in\mod B$ whose support is in $Q_0\setminus F_0$, the object $M\lten_B T$ is zero in $\Dd^{\rm b}(\Ee)$. Consequently the object $\A\lten_B T$ is zero in $\Dd(\A^{\rm op}\otimes\Ee)$, and the morphism $$\xymatrix{ \RHom_A(DA,\bar{A})\lten_A T[2]\ar@<-.5ex>[r] & \A\lten_A
T}$$ is an isomorphism in $\Dd^{\rm b}(\A^{\rm op}\otimes\Ee)$.

\end{proof}
\begin{lema}\label{lem2}
There exists a morphism \[ \xymatrix{\alpha:\RHom_{A}(DA,\A)\lten_A T\ar@<-.5ex>[r] & \RHom_{A}(DA,\A)\lten_A (\A\lten_A T)}\quad\textrm{ in }\Dd^{\rm b}(\A^{\rm op}\ten \Ee)\] whose cone is in $\Dd^{\rm b}(\A^{\rm op}\ten \Pp)$.
\end{lema}
\begin{proof}
Applying the functor $\RHom_A(DA,\A)\lten_A -\lten_A T$ to the exact sequence $$\xymatrix{Ae_FA\; \ar@{>->}[r] & A\ar@{->>}[r] & A/Ae_F A=\A }$$ of $A$-bimodules, we get the following triangle 
$$\xymatrix{ \RHom_A(DA,\A)\lten_A e_F T\ar@<-.5ex>[r]& \RHom_A(DA,\A)\lten_A T\ar@<-.5ex>[r] &
\RHom_{A}(DA,\A)\lten_A (\A\lten_A T)\ar@<-.5ex>[r] & }.$$
Since $A$ is finite dimensional and of finite global dimension, $\RHom_A(DA,\A)$ is in ${\sf thick}(A)$, the thick subcategory of $\Dd(A)$ generated by $A$. Thus $\RHom_A(DA,\A)\lten_A e_F T$ is in ${\sf thick} (e_F T)$, the thick subcategory of $\Dd(\Ee)$ generated by $e_F T$. By (H1), $e_FT$ is projective injective, and therefore it is in $\Pp$. Hence the cone of the morphism \[ \xymatrix{\RHom_{A}(DA,\A)\lten_A T\ar@<-.5ex>[r] & \RHom_{A}(DA,\A)\lten_A (\A\lten_A T)}\]  is in $\Dd^{\rm b}(\A^{\rm op}\ten \Pp)$.
\end{proof}
\begin{lema}\label{lem3}
There is an isomorphism \[ \xymatrix{\RHom_{A}(DA,\A)\lten_A \A \ar@<-.5ex>[r]^-\sim & \RHom_{\A}(D\A,\A)} \quad \textrm{ in } \Dd^{\rm b}(\A^{\rm op}\ten A).\]
\end{lema}

\begin{proof}
Let $Z$ be an injective resolution of $\A$ as a right $A$-module. Then we have $$\RHom_{A}(DA,\A)\lten_A \A\simeq \Hom_A(DA,Z)\otimes_A\A$$ using the fact that $\Hom_A(DA,Z)$ is a complex of projectives.

For each $A$-module $M$, we define a morphism
$H_M:\Hom_A(DA,M)\otimes_A\A\rightarrow \Hom_A(D\A,M)$ as follows.
Let $\varphi$ be in $\Hom_A(DA,M)$, $\overline{a}\in\A$, and $\Phi\in D\A$. We define $H_M(\varphi\otimes \overline{a})(\Phi):=\varphi (\Phi_{\overline{a}})\in M,$ where $\Phi_{\overline{a}}(b)=\Phi(b\overline{a})$ for any $b\in A$. The morphism $H_M$ is functorial in $M$. An easy computation shows that the map $H_{DA}$ is an isomorphism. Since $Z$ is a bounded complex of injective $A$-modules, we get the following isomorphism in $\Dd^{\rm b}(\A^{\rm op}\otimes A)$
$$\xymatrix{\RHom_{A}(DA,\A)\lten_A \A= \Hom_A(DA,Z)\otimes_A\A \ar[rr]^(.55){H_Z}& &\Hom_A(D\A,Z)=\RHom_A(D\A,\A).}$$
Moreover by Lemma~\ref{restriction} any injective resolution of an $\A$-module $X$ is an injective resolution of $X$ viewed as an $A$-module, and the restriction functor is fully faithful. Thus we get isomorphisms in  $\Dd^{\rm b}(\A^{\rm op}\otimes A)$
$$\RHom_A(D\A,\A)=\Hom_A(D\A,Z)\simeq\Hom_{\A}(D\A,Z)=\RHom_{\A} (D\A,\A).$$
\end{proof}

\begin{proof}[Proof of Proposition~\ref{construction alpha}]
Combining Lemmas~\ref{lem1},~\ref{lem2} and~\ref{lem3}, we get the following morphisms in $\Dd^{\rm b}(\A^{\rm op}\ten \Ee)$
\[ \xymatrix{ \A\lten_A T\ar@{-}_\wr[d]^-{\rm Lemma~\ref{lem1}} &&  \\
 \RHom_A(DA,\A)\lten_A T[2] \ar@<-.5ex>[rr]_-{\rm Lemma~\ref{lem2}}^-\alpha && \RHom_A(DA,\A)\lten_A(\A\lten_A T)[2]\ar@{-}_\wr[d]^-{\rm Lemma~\ref{lem3}}  \\
&&\RHom_{\A}(D\A,\A)\lten_A T[2]\ar@{-}_\wr[d]\\  && \RHom_{\A}(D\A,\A)\lten_{\A}(\A\lten_A T)[2]
}\]
By Lemma~\ref{lem2},  the cone of the horizontal morphism is in $\Dd^{\rm b}(\A^{\rm op}\ten\Pp)$. Hence we get Proposition~\ref{construction alpha}.
\end{proof}

Using Proposition~\ref{construction alpha}, we can apply Proposition~\ref{univprop} for $\Lambda=\A$ and $M=\A\lten_A T$ and we obtain the following.

\begin{cora}\label{existence G}
In the setup of Theorem~\ref{theoremgrading}, there exists a triangle functor $G:\Cc_{\A}\rightarrow \underline{\Ee}$ such that the following diagram commute:$$\xymatrix{\Dd^{\rm b}(\A)\ar[r]_{R}\ar[d]^{\pi_{\A}}\ar@/^.5cm/[rrr]^{-\lten_{\A}(\A\lten_A T)} & \Dd^{\rm b}(A)\ar[rr]_{-\lten_A
    T}& &
  \Dd^{\rm b}(\Ee) \ar[d]  \\ \Cc_{\A}\ar[rrr]^{G}&&&\underline{\Ee}.}$$
\end{cora}

\subsection{Proof of the equivalence}

In this subsection we show that the triangle functor $G:\Cc_{\A}\rightarrow \underline{\Ee}$ of Corollary~\ref{existence G}  is an equivalence.  The proof is separated into four steps.  First we show that the functor $G$ sends the cluster-tilting object $\A$ to the cluster-tilting object $T\in\underline{\Ee}$. The second step (Proposition~\ref{prop injective map}) consists of proving that $G$ induces an injective map from the endomorphism algebra of $\A\in\Cc_{\A}$ to the endomorphism algebra of $T\in\underline{\Ee}$. In Proposition~\ref{algebraBbar}, we prove that these two endomorphism algebras are isomorphic. The last step follows from \cite{Kel4} (Proposition~\ref{keller reiten}) and finishes the proof of Theorem~\ref{theoremgrading}.
\medskip

We have a triangle 
$$\xymatrix{Ae_F A\lten_A T\ar@<-.5ex>[r] & T\ar@<-.5ex>[r] & \A\lten_A T\ar@<-.5ex>[r] & Ae_F A\lten_A T[1]}\quad \textrm{ in } \Dd(\A^{\rm op}\ten \Ee)$$  and the object
 $Ae_FA\lten_A T=e_FT$ is in $\Pp$. Hence $\A\lten_A T$  is isomorphic to $T$ in $\underline{\Ee}$. Therefore
the triangle functor $G$ sends the cluster-tilting object $\A$ to the cluster-tilting object $T$.

\medskip

In order to prove Proposition~\ref{prop injective map} which is the second step, we first prove the following.

\begin{lema}\label{lema injectivity}
Let $M$ be an object in $\Dd^{\rm b}(\A)$ such that its homology is concentrated in non positive degrees, and let $f\in\Hom_{\Dd^{\rm b}(\A)}(\A,M)$ be a non zero morphism. Then the morphism $f\lten_A T\in \Hom_{\Dd^{\rm b}(\Ee)}(\A\lten_A T,M\lten_A T)$ is non zero viewed as a morphism in $\underline{\Ee}$.
\end{lema}

\begin{proof}
By Lemma~\ref{restriction} $(c)$, the restriction functor $\Dd^{\rm b}(\A)\rightarrow \Dd^{\rm b}(A)$ is fully faithful. We first show that the morphism $f\lten_A T$ is non zero as a morphism in $\Dd^{\rm b}(\Ee)$.

Now by Lemma~\ref{restriction} $(a)$, $\A$ is quasi-isomorphic to some complex in $\Dd^{\rm b}(A)$ of the form:
\[ \xymatrix{\A= \cdots\ar[r] & 0\ar[r] & Q^{-2}\ar[r] & Q^{-1}\ar[r] & (1-e_F) A\ar[r] & 0 \ar[r] & \cdots}\] with $Q^{-1}$ and $Q^{-2}$ in $\add (e_F A)$.
Since $M$ is concentrated in non positive degrees, $M$ is quasi isomorphic to some bounded complex in $\Dd^{\rm b}(A)$ of the form
\[ \xymatrix{M=\cdots\ar[r] & 0\ar[r] &\cdots \ar[r] & e_{-3}A\ar[r] & e_{-2}A\ar[r] & e_{-1} A\ar[r] & e_0A\ar[r] & 0\ar[r] &\cdots}\]
Hence $f$ is isomorphic to a morphism of complexes of the form
\[\xymatrix{ \cdots\ar[r] & 0\ar[r]\ar[d] & Q^{-2}\ar[r]\ar[d]^{f^{-2}} & Q^{-1}\ar[d]^{f^{-1}}\ar[r] & (1-e_F) A\ar[r]\ar[d]^{f^0} & 0 \ar[d]\ar[r] & \cdots\\ 
\cdots\ar[r] & e_{-3}A\ar[r]^{d^{-3}} & e_{-2}A\ar[r]^{d^{-2}} & e_{-1} A\ar[r]^{d^{-1}} & e_0A\ar[r] & 0\ar[r] &\cdots}.\]
Hence $f°$ is isomorphic to $f\circ p$,  where $p:(1-e_F)A\rightarrow \A$ is the projective cover. Hence $f^0$ is a non zero morphism.

 By definition the morphism $f^0\lten_A T$ is the morphism of complexes:
\[\xymatrix{ \cdots\ar[r] & 0\ar[r]\ar[d] & 0\ar[r]\ar[d] & 0\ar[d]\ar[r] & (1-e_F)T\ar[r]\ar[d]^{f^0\ten_A T} & 0 \ar[d]\ar[r] & \cdots\\ 
\cdots\ar[r] & e_{-3}T\ar[r]^{d^{-3}} & e_{-2}T\ar[r]^{d^{-2}} & e_{-1} T\ar[r]^{d^{-1}} & e_0T\ar[r] & 0\ar[r] &\cdots}\]
Note that $f^0\lten_A T=f^0\ten_A T$ since $(1-e_F)A$ and $e_0A$ are projective $A$-modules.
In this diagram, the maps $d^i$ and $f^0\ten_A T$ are morphisms of degree 0, since they come from morphisms in $\add(A)$.

Denote by $\xymatrix{P:=\cdots\ar[r] & P^{-2}\ar[r] & P^{-1}\ar[r] & P^0}$ a projective resolution in $\Ee$ of the object $(1-e_F)T$ and denote by $p'$ the map $p': P\rightarrow (1-e_F)T$. Note that the $P^i$ are also injective since $\Ee$ is a Frobenius category. Assume that the morphism $f^0\lten_A T$ vanishes in $\Dd^b(\Ee)$. It implies that the morphism $(f^0\ten_A T)\circ p'$ is homotopic to zero. 
\[\xymatrix{ \cdots\ar[r] & P^{-3}\ar[r]\ar[d] & P^{-2}\ar[ddl]\ar[r]\ar[d] & P^{-1}\ar[r]\ar[d] \ar[ddl]& P^0\ar[r] \ar[d]^{p'} \ar[ddl]& 0\ar[r]\ar[d] & \cdots\\
\cdots\ar[r] & 0\ar[r]\ar[d] & 0\ar[r]\ar[d] & 0\ar[d]\ar[r] & (1-e_F)T\ar[r]\ar[d]^{f^0\ten_A T}\ar@{..>}[dl]^x & 0 \ar[d]\ar[r] & \cdots\\ 
\cdots\ar[r] & e_{-3}T\ar[r]^{d^{-3}} & e_{-2}T\ar[r]^{d^{-2}} & e_{-1} T\ar[r]^{d^{-1}} & e_0T\ar[r] & 0\ar[r] &\cdots}\]
Since the complex $M\lten_A T$ is quasi isomorphic to a bounded complex with components in $\add(T)$, an easy induction shows that the map $f^0\ten_A T=d^{-1}x$ is homotopic to zero. Since $f^0\ten_A T$ and $d^{-1}$ are homogeneous of degree 0, then we can assume that  $x:(1-e_F)T\rightarrow e_{-1}T$ is homogenous of degree 0. Therefore $x$ can be written $y\ten_A T$ for a map $y\in\Hom_{\Dd^b(A)}((1-e_F)A,e_0A)$. Therefore $f^0=d^{-1}y$ is homotopic to zero, hence $f^0$ vanishes as a map in $\Dd^b(A)$. Consequently,  if $f^0$ is not zero in $\Dd^b A$, then the map $f\lten_A T$ is not  zero in $\Dd^b(\Ee)$, since $f^0\lten_A T=(f\lten_A T)\circ(p\lten_A T)$.

Assume now that $f\lten_A T$ is zero in $\underline{\Ee}$. Since $p\lten_A T: (1-e_F)T\rightarrow \A\lten_A T$ is an isomorphism in $\underline{\Ee}$, then $f\lten_A T$ vanishes in $\underline{\Ee}$ if and only if $f^0\lten_A T$ vanishes in $\underline{\Ee}$. Thus $f^0\lten_A T$ factors through an object $P$ in $\Dd^{\rm b}(\Pp)$. 
For $i\in\mathbb{Z}$ we denote by $P^{\geq i}$ and $P^{\leq i}$  the positive and negative truncations
$$\xymatrix@-.2cm{P^{\geq i}:=\cdots 0\ar[r] & P^i\ar[r] & P^{i+1}\ar[r] & P^{i+2}\ar[r] & \ldots }\quad \xymatrix@-.2cm{P^{\leq i}:=\cdots \ar[r] &P^{i-2}\ar[r] & P^{i-1}\ar[r] & P^i\ar[r] & 0\cdots}.$$  For any $i\in\mathbb{Z}$ the object $P^i$ is projective and injective, thus the complex $P$ and all its truncations are fibrant and cofibrant. Therefore the space $\Hom_{\Dd^{\rm b}(\Ee)}((1-e_F)T,P^{\leq -1})$ vanishes because $(1-e_F)T$ is concentrated in degree $0$. And the space $\Hom_{\Dd^{\rm b}(\Ee)}(P^{\geq 1},M\lten_A T)$ vanishes since $M$ (hence $M\lten_A T$) is concentrated in non positive degrees. Consequently we can assume that $P=P^0$ is a stalk complex with $P\in\add(e_FT)=\Pp$. We write $f^0\lten_A T=g\circ h$ with $g\in\Hom_{\Dd^{\rm b}(\Ee)}(P^0,e_0 T)$ and $h\in\Hom_{\Dd^{\rm b}(\Ee)}((1-e_F)T,P^0)$. Now by definition $f^0\lten_A T\in\Hom_{\Dd^{\rm b}(\Ee)}((1-e_F)T,e_0T)$ is of degree 0. By hypothesis $(H4)$, the morphism $h$ is  in non negative degrees. Since $B=\End_{\Dd^{\rm b}(\Ee)}(T)$ is only in positive degree we get a contradiction.

\end{proof} 

Using this fundamental lemma, we can prove the following.
\begin{prop}\label{prop injective map}
The functor $G$ constructed in subsection 3.3 induces an injective map
\[ \xymatrix{\End_{\Cc_{\A}}(\A)\ar[r] & \underline{\End}_{\Ee}(T).} \]
\end{prop}

\begin{proof}
The two algebras $\End_{\Cc_{\A}}(\A)=\bigoplus_{i\geq 0}\Hom_{\Dd^{\rm b}(\A)}(\A,\mathbb{S}^{-i}\A[2i])$ and $\bar{B}=\underline{\End}_{\Ee}(T)$ are graded algebras. The first part of the proof consists of showing that the functor $G:\Dd^b(\A)\rightarrow \underline{\Ee}$, which sends $\A\in\Cc_{\A}$ to $T\in\underline{\Ee}$, induces a morphism of graded algebras $\xymatrix{\End_{\Cc_{\A}}(\A)\ar[r] & \underline{\End}_{\Ee}(T)} $ and then that this morphism is injective using Lemma~\ref{lema injectivity}.

\medskip

Let $i\geq 0$ and $f^i\in\Hom_{\Dd^{\rm b}(\A)}(\A,\mathbb{S}^{-i}\A[2i])$ be a non zero morphism.  We define $\bar{\theta}^p$ for $p\geq 0$ by induction by $\bar{\theta}^1=\bar{\theta}=\RHom_{\A}(D\A,\A)[2]=\mathbb{S}_2^{-1}A[2]$ and $\bar{\theta}^p:=\bar{\theta}\lten_{\A}\bar{\theta}^{p-1}$. Hence $\bar{\theta}^p=\mathbb{S}^{-p}(\A)[2p]$ and $f^i$ is a non zero map in $\Hom_{\Dd^{\rm b}(\A)}(\A,\bar{\theta}^i)$. 

We first show that the triangle of Proposition~\ref{triangle} is in fact a triangle 
\[\xymatrix{ \RHom_A (DA,\A)\lten_A B(-1)[2]\ar[r] & \A\lten_A B\ar[r] & \A \ar[r] &  \RHom_A (DA,\A)\lten_A B(-1)[3]}\]
in the category $\Dd^b(\gr (\A\ten B))$, where $(1)$ is the degree shift of the graded algebra $B$. Indeed the map $f:X\rightarrow Y$ constructed in subsection~2.4 is homogeneous of degree 1. Therefore by Lemma~\ref{lem1}, we obtain an isomorphism 
$$\A\lten_A B (1)\lten_B T\simeq \RHom_A (DA,\A)\lten_A T[2].$$
Since the maps constructed in Lemmas~\ref{lem2} and~\ref{lem3} come from maps in $\Dd^b(A)$ they are of degree $0$.  Consequently  the triangle of Proposition~\ref{triangle} is a triangle\[\xymatrix{ \RHom_A (DA,\A)\lten_A B(-1)[2]\ar[r] & \A\lten_A B\ar[r] & \A \ar[r] &  \RHom_A (DA,\A)\lten_A B(-1)[3]} \] in $\Dd^b(\gr (\A\ten B))$. Hence  by  Proposition~\ref{construction alpha} we get a map $\alpha:\A\lten_AB(1)\lten_B T\rightarrow \bar{\theta}\lten_A T$ in $\Dd^{\rm b}(\A^{\rm op}\ten\Ee)$ whose cone is in $\Dd^{\rm b}(\Pp)$. We denote by $\alpha^i$ the composition:\[ \xymatrix{\alpha^i: \A\lten_A B(i)\lten_BT\ar[r]^-{\alpha} & \bar{\theta}\lten_A B(i-1)\lten_B T\ar[rr]^{\bar{\theta}^0\lten_{\A}\alpha} && \bar{\theta}^2\lten_A B(i-2)\lten_B T \ar[r] & \cdots \ar[r] & \bar{\theta}^i\lten_{A}T.}\] 
The image of $f^i\lten_A T$ in $\underline{\Ee}$ is isomorphic (as a morphism in $\underline{\Ee}$) to the left fraction $(\alpha^i)^{-1}\circ (f^i\lten_A T)$:
\[\xymatrix{\A\lten_A T\ar[dr]_{f^i\lten_A T} && \A\lten_AB(i)\lten_B T\ar[dl]^{\alpha^i}\\ & \bar{\theta}^i\lten_A T&}\]
Thus the image of $f^i\lten_A T$ in $\underline{\Ee}$ is a map of degree $i$ and $G$ induces a morphism of graded algebras $\xymatrix{\End_{\Cc_{\A}}(\A)\ar[r] & \underline{\End}_{\Ee}(T)} $

\medskip

The functor $-\lten_A T:\Dd^b\A\rightarrow \underline{\Ee}$ induces a map $\Hom_{\Dd^{\rm b}(\A)}(\A,\bar{\theta}^i)\rightarrow \bar{B}$. Since the object $\bar{\theta}^i\in\Dd^{\rm b}(\A)$ has its homology concentrated in non positive degrees, this map is injective by  Lemma~\ref{lema injectivity}. By the remark above, if $f^i\in\Hom_{\Dd^b\A}(\A,\bar{\theta}^i)$, then $f^i\lten_AT$ is a morphism of degree $i$ in the graded algebra $\underline{\End}_{\Ee}(T)=\bar{B}$. Hence there is an injective map $\Hom_{\Dd^{\rm b}(\A)}(\A,\bar{\theta}^i)\rightarrow \bar{B}_i$, where $\bar{B}_i$ is the degree $i$ part of the graded algebra $\bar{B}=\underline{\End}_{\Ee}(T)$.  Therefore the functor $G:\Cc_{\A}\rightarrow \underline{\Ee}$ induces a morphism of graded algebras \[\bigoplus_{i\geq 0}\Hom_{\Dd^{\rm b}(\A)}(\A,\mathbb{S}^{-i}\A[2i])=\bigoplus_{i\geq 0}\Hom_{\Dd^{\rm b}(\A)}(\A,\bar{\theta}^i)\longrightarrow \bar{B}=\bigoplus_{i\geq 0}\bar{B}_i,\] which is an injection $\Hom_{\Dd^{\rm b}(\A)}(\A,\bar{\theta}^i)\rightarrow \bar{B}_i$ for any $i\geq 0$. Consequently it is an injective algebra morphism.

\end{proof}

The next result will be used in the proof of Proposition~\ref{algebraBbar} which is the third step.
\begin{lema}
Let $\Jac(Q,W,F)$ be a frozen Jacobian algebra with $W$ reduced. Assume there is a grading $\varphi:Q_1\rightarrow \{0,1\}$ satisfying the hypotheses $(H2)$, $(H3)$ and $(H4)$.   We denote by $\bar{Q}$ the full subquiver of $Q$ with set of vertices $\bar{Q}_0:=Q_0\setminus F_0$, and by $\bar{W}$ the image of $W$ under the projection $kQ\rightarrow k\bar{Q}$. Then the set $\{\partial_a\bar{W}, \varphi(a)=1\}$ is linearly independent. In particular, $\partial_a\bar{W}$ does not vanish for $a$ in $\bar{Q}_1$ with $\varphi(a)=1$.
\end{lema}

\begin{proof}
Let $a\in \bar{Q}_1$ be an arrow with $\varphi(a)=1$. By condition (H3) the potential $\bar{W}$ is homogeneous of degree~1, and by assumption the degree map $\varphi$ has non-negative
 values. Thus any term in the potential $\bar{W}$ contains exactly one arrow of degree 1.  Consequently a cycle containing $a$ which is a summand of $W$ does not contain any other arrow of degree 1. Then by (H4) it follows that this cycle does not pass through vertices in $F_0$. Therefore we have $\partial_aW=\partial_a \bar{W}$. Now by (H2) the set  $\{\partial_a W, \varphi(a)=1\}$ is linearly independent, therefore we get the result.
\end{proof}

\begin{prop}\label{algebraBbar}
Under the assumptions (H1)-(H4) there is an isomorphism of algebras
$$ \End_{\underline{\Ee}}(T)\simeq \End_{\Cc_{\A}}(\A).$$
\end{prop}
\begin{proof}
The algebra $\bar{B}=B/Be_FB$ is isomorphic to the Jacobian algebra $\Jac(\bar{Q},\bar{W})$ (cf section~2), where $\bar{Q}$ is the full subquiver of $Q$ whose vertices are not in $F_0$ and where $\bar{W}$ is the image of $W$ under the projection $\xymatrix{kQ\ar@{->>}[r] & kQ/kQe_FkQ\simeq k\bar{Q}}.$  

Let $Q'$ be the subquiver of $Q$ defined by:
\begin{itemize}
\item $Q'_0:=Q_0$;
\item $Q'_1:=\{a\in Q_1, \varphi(a)=0\}$.
\end{itemize}

By definition we have $A\simeq kQ'/\langle \partial_a W, a\notin F_1
\textrm{ and }
\varphi(a)=1\rangle.$
Let $\bar{Q}'$ be the full subquiver of $Q'$ with vertices $\bar{Q}'_0=Q'_0\setminus F_0$. Thus we get $\A\simeq k\bar{Q'}/\langle \partial_a \bar{W},
\varphi(a)=1\rangle.$ 

By Theorem~\ref{keller}, 
the endomorphism algebra
$\End_{\Cc_{\A}}(\A)$ is isomorphic to the Jacobian algebra $\Jac(\widetilde{\bar{Q}'},W_{\A})$ where $\widetilde{\bar{Q}'}$ and $W_{\A}$ are defined in Theorem~\ref{keller}.   By the previous lemma, the set $\{\partial_a\bar{W}, \varphi(a)=1\} $ is a basis for the ideal of relations of $\A$. Hence we immediately see that $\widetilde{\bar{Q}'}=\bar{Q}$ and that $W_{\A}$ is the potential  $$W_{\A}=\sum_{a,\varphi(a)=1} a\partial_a \bar{W}.$$
Therefore the potential $\bar{W}$ is cyclically equivalent to the potential $W_{\A}$. Thus we have an isomorphism
$\Jac(\widetilde{\bar{Q}'},W_{\A})\simeq \Jac(\bar{Q},\bar{W})$
which gives the desired isomorphism.  

\end{proof}

Now we finish the proof of Theorem~\ref{theoremgrading}. The triangle functor $G:\Cc_{\A}\rightarrow \underline{\Ee}$ constructed in subsection 3.3 sends the cluster-tilting object $\A\in\Cc_{\A}$ to the cluster-tilting object $T\in\underline{\Ee}$. By Proposition~\ref{prop injective map}, it induces an injective map $\End_{\Cc_{\A}}(\A)\rightarrow \underline{\End}_{\Ee}(T)$. These two algebras are finite dimensional algebras which are isomorphic by Proposition~\ref{algebraBbar}. Therefore $G$ induces a bijection $\End_{\Cc_{\A}}(\A)\simeq \underline{\End}_{\Ee}(T)$. Then we conclude the proof of Theorem~\ref{theoremgrading} using the following.

\begin{prop}[\cite{Kel4} Lemma 4.5]\label{keller reiten}
Let $\Cc$ and $\Cc'$ be $\Hom$-finite $2$-Calabi-Yau triangulated categories. Let
$T$ (resp. $T'$) be a cluster-tilting object in $\Cc$
(resp. $\Cc'$). If we have a triangle functor $G:\Cc\rightarrow \Cc'$
which sends $T$ to $T'$ and which induces an isomorphism between
$\End_\Cc(T)$ and $\End_{\Cc'}(T')$, then $G$ is an equivalence.
\end{prop}

\section{2-Calabi-Yau categories associated with elements in the Coxeter group}

In this section we apply Theorem~\ref{theoremgrading} to the categories associated with elements in the Coxeter group introduced in \cite{Bua2}. 

\subsection{Results of \cite{Bua2} and \cite{BIRSm}}
Let $Q$ be a finite quiver without oriented cycles. We denote  as usual by $Q_0=\{1,\ldots, n\}$ the set of vertices and by $Q_1$ the set of arrows. The preprojective algebra associated to $Q$ is the algebra
$$k\overline{Q}/\langle \sum_{a\in Q_1} aa^*-a^*a\rangle$$
where $\overline{Q}$ is the double quiver of $Q$, which by definition is obtained from $Q$ by adding to each arrow 
$a:i\rightarrow j$ in $Q_1$ an arrow $a^*:i\leftarrow j$ pointing in the opposite direction. 
We denote by $\Lambda$ the completion of the preprojective algebra associated to $Q$, and by $\fl\Lambda$ 
the category of right $\Lambda$-modules of finite length. 

Let $C_Q$ be the 
Coxeter group associated to $Q$. It is defined by the generators $s_i$ where $i\in Q_0$ and by the relations
\begin{itemize}
 \item $s_i^2=1$,
\item $s_is_j=s_js_i$ if there is no arrow between $i$ and $j$,
\item $s_is_js_i=s_js_is_j$ if there is exactly one arrow between $i$ and $j$.
\end{itemize}

 A \emph{reduced expression} $\ww=s_{u_1}\ldots s_{u_l}$ of an element $w$ in $C_Q$ is an expression with $l$
 smallest possible. When  $\ww=s_{u_1}\ldots s_{u_l}$ is a reduced expression of $w$, the integer $l(w):=l$ is called the length of $w$. 

For a vertex $i$ in $Q_0$ we denote by $\Ii_i$ the ideal $\Lambda(1-e_i)\Lambda$. Let $\ww=s_{u_1}\ldots s_{u_l}$ be a reduced expression of an element $w$ in $C_Q$. For $p\leq l$ we denote by 
$\Ii_{\ww_p}$ the ideal $\Ii_{u_p}\Ii_{u_{p-1}}\ldots\Ii_{u_1}$. The ideal $\Ii_{\ww}:=\Ii_{\ww_l}$ depends only on the element $w\in C_Q$ and not on the choice of the reduced expression $\ww$. Therefore we denote by $\Lambda_{w}$ the algebra 
$\Lambda/\Ii_{\ww}$ and by $\Ee_{w}:=\Sub\Lambda_{w}$ the subcategory of $\fl\Lambda$ consisting of submodules of finite direct sums of copies of $\Lambda_w$.

We have the following \cite[Theorem III.2.8]{Bua2}.
 
\begin{thma}[Buan-Iyama-Reiten-Scott]\label{birsc}
 Let $w$ be an element in the Coxeter group $C_Q$. Then we have the following.
\begin{itemize}
\item[(a)]The category $\Ee_w$ is a $\Hom$-finite Frobenius stably 2-CY category. 
\item[(b)]
For any reduced expression $\ww=s_{u_1}\ldots s_{u_l}$ of $w$, the object $T_\ww=\bigoplus_{p=1}^le_{u_p}(\Lambda/\Ii_{\ww_p})$ 
is a cluster-tilting object. 
\item[(c)] The projective-injective indecomposable objects are $e_{u_{t_i}}(\Lambda/\Ii_{\ww_{t_i}})$  where $t_i$ is the maximal integer such that $u_{t_i}=i$ for $i\in Q_0$.
\end{itemize}
\end{thma}
The cluster-tilting object $T_\ww$ depends on the choice of the reduced expression of $w$. 
We refer to a cluster-tilting object of this form as a \emph{standard cluster-tilting object}. 
Note that by mutation we may get other cluster-tilting objects which are not standard.

We now define a quiver $Q_\ww$ associated with a reduced expression $\ww=s_{u_1}\ldots s_{u_l}$ of an element $w\in C_Q$ as follows:

\begin{itemize}
\item vertices: $1,\ldots,l(w)$.
\item for each $i\in Q_0$, one arrow $t\leftarrow s$ if $t$ and $s$ are two consecutive vertices of type $i$ 
(\emph{i.e.} $u_s=u_t=i$) and $t<s$ (we call these arrows \emph{arrows going to the left});
\item for each $a:i\rightarrow j \in Q_1$, put $a:t\rightarrow s$ if $t$ is a vertex of type $i$, 
$s$ of type $j$, and if 
  there is no vertex of type $i$ between $t$ and $s$, and if $s$ is the 
last vertex of type $j$ before the next vertex of type $i$ in the expression
  $\ww=s_{u_1}\ldots s_{u_l}$ (we call these arrows the \emph{$Q$-arrows})
\item for each $a:i\rightarrow j\in Q_1$, put $a^*: t\rightarrow
  s$ if $t$ is of type $j$, $s$ is of type $i$, if there is no vertex of type $j$ between $t$ 
and $s$
  and if $s$ is the last vertex of type $i$ before the next vertex of type $j$ in the expression $\ww$ (we call these arrows the \emph{$Q^*$-arrows}).
\end{itemize}

For each $Q$-arrow $a:t\rightarrow s$ in $Q_\ww$,
 we denote by $W_{a}$ the
composition $aa^*p$ if there is a (unique) $Q^*$-arrow
$a^*:r\rightarrow t$ in $Q_\ww$ where $u_r=u_s$ and where $p$ is the
composition of arrows going to the left $r\leftarrow \cdots \leftarrow s$. Otherwise we put
$W_a=0$. For
each $Q^*$-arrow $a^*:t\rightarrow r$ in $Q_\ww$, we denote by $W_{a^*}$ the
composition $a^*a p$ if there exists a (unique) $Q$-arrow
$a:s\rightarrow t$ with $u_s=u_r$ in $Q_\ww$ and where $p$ is the composition of arrows going to the left
$s\leftarrow \cdots \leftarrow r$. Otherwise we put
$W_{a^*}=0$. Then let $W_\ww$ be the sum
$$W_\ww=\sum_{a\ Q\textrm{-arrow}}W_{a}-\sum_{a^*\ Q^*\textrm{-arrow}}W_{a^*}.$$

We have the following~\cite[Theorem 6.8]{BIRSm}.

\begin{thma}[Buan-Iyama-Reiten-Smith]\label{birs}
Let $\ww=s_{u_1}\ldots s_{u_l}$ be a reduced expression of an element $w$ of the Coxeter group $C_Q$. Let $T_\ww$ be the associated standard cluster-tilting object in the category $\Ee_w$. Then we have an isomorphism $$\End_{\Ee_w}(T_\ww)\simeq Jac(Q_\ww,W_\ww,F)$$ where $F_0:=\{t_1,\ldots, t_n\}$ consists of the last vertex of each kind in the reduced expression $\ww$ and 
$F_1:=\{a\in Q_1, s(a)\in F_0 \textrm{ and }t(a)\in F_0\}$.
\end{thma}

\subsection{Definition of the grading}

The algebra $\Lambda$ and the associated Coxeter group do not depend on the orientation of $Q$. For any reduced expression $\ww=s_{u_1}\ldots s_{u_l}$ of $w\in C_Q$, the category $\Ee_w$ and the cluster-tilting object $T_\ww$ do not depend on the orientation of $Q$. From now on we assume that the orientation of $Q$ satisfies the property that $t_i<t_j$ if there exists an arrow $i\rightarrow j$, where $t_i$ is the maximal integer satisfying $u_{t_i}=i$.

We define a grading on the quiver $Q_\ww$ as follows:
\begin{itemize}
\item $\varphi(b)=1$ if $b$ is a $Q^*$-arrow;
\item $\varphi(b)=0$ if $b$ is a $Q$-arrow or an arrow going to the left.
\end{itemize}

With this choice of grading we show that our axioms are satisfied.
\begin{lema}
The graded Jacobian algebra $(\Jac(Q_\ww,W_\ww,F),\varphi)$ satisfies the conditions $(H1)$-$(H4)$ of Theorem~\ref{theoremgrading}.
 \end{lema}
\begin{proof}
(H1) This holds by Theorems~\ref{birsc} and~\ref{birs}. 

\smallskip
\hspace{.5cm}(H2) The potential $W_\ww=W$ is reduced. Moreover it is easy to see that two different terms of the
potential $W$ differ by at least two arrows. It follows that  the set $\{\partial_aW|\ a\notin F_1, \partial_aW\neq 0\}$ is a basis for the ideal of relations. We next show that for any arrow $a\notin F_1$, the derivative $\partial_aW$ does not vanish. Assume that some $Q$-arrow or $Q^*$-arrow $a:r\rightarrow s$ does not appear in the potential, where $u_r=i$, $u_s=j$ and there is an arrow between $i$ and $j$ in the quiver $Q$. Then  there is no $u_t$ of type $i$ with $t>r$. Thus we have $r=t_i$. Then we must have $s=t_j$. Therefore $a:r=t_i\rightarrow s=t_j$ is in $F_1$. Now let $p:r\leftarrow s$ be an arrow going to the left, where $r$ and $s$ are two consecutive vertices of type $i$. Since the expression $s_{u_1}\ldots s_{u_l}$ is reduced, there exists $t$ with $r<t<s$ such that $u_t$ is of type $j$ and there is (at least) one arrow between $i$ and $j$ in the quiver $Q$. Let $t$ be the maximal integer with this property. Then there is an arrow $a:r\rightarrow t$, and there is also an arrow $a^*:t\rightarrow v$ where $v\geq s$ and $u_v$ is of type $i$. Thus $p$ appears in the potential $W$. Therefore for any arrow $a$ of $Q_\ww$ which is not in $F_1$, the derivative $\partial_aW$ is not zero, and we have (H2).       

\smallskip
\hspace{.5cm}(H3) By the definition of the potential in Theorem~\ref{birs}, this follows immediately.

\smallskip
\hspace{.5cm}(H4) By the choice of the orientation of $Q$, any arrow in $F_1$ is a $Q$-arrow, hence of degree 0.
 Thus all arrows with target in $F_0$ and source not in $F_0$ are $Q^*$-arrows, hence of degree $1$, and all arrows with source in $F_0$ and target not in $F_0$ are arrows going to the left, hence of degree $0$. This implies condition $(H4)$.

\end{proof}

The following direct consequence is one of our main results.

\begin{thma}\label{applicationbirs}
The stable category $\underline{\Ee}_w$ is a generalized cluster category.
\end{thma}

\subsection{Meaning of the grading}
We show that the grading $\varphi$ on $\End_{\Ee_w}(T)$ defined in the previous subsection is induced by a natural grading on the preprojective algebra. 

Let $\ww=s_{u_1}\ldots s_{u_l}$ be a reduced expression of an element $w$ of the Coxeter group $C_Q$. Assume that the orientation of $Q$ satisfies the property that $t_i<t_j$ if there exists an arrow $i\rightarrow j$,
where $t_i$ is the maximal integer satisfying $u_{t_i}=i$. 

We define a grading on the double quiver $\overline{Q}$ as follows:
\begin{itemize}
\item $deg(a)=0$ if $a$ is an arrow of $Q$;
\item $deg(a^*)=1$ if $a^*:j\rightarrow i$ is an arrow pointing in the opposite direction of an arrow $a:i\rightarrow j$ in $Q$.
\end{itemize}
The ideal of relations $(\sum_{a\in Q_1}aa^*-a^*a)$ is
homogeneous of degree $1$, thus the grading on the double quiver $\overline{Q}$ induces a grading on the preprojective algebra $\Lambda$.

For any $i$ in $Q_0$, the $\Lambda$-module $e_i\Lambda$ can be seen as a graded $\Lambda$-module with top in degree $0$. 
Then the ideal $\Ii_i=\Lambda(1-e_i)\Lambda$ is a
graded ideal.  
Hence for $p\leq l$ the ideal $\Ii_{\ww_p}=\Ii_{u_p}\ldots \Ii_{u_1}$ is a graded
ideal and the module $T_p=e_{u_p}(\Lambda/\Ii_{\ww_p})$ is a finite
length graded $\Lambda$-module.
Therefore the standard cluster-tilting object $T_\ww=T_1\oplus\cdots\oplus T_l$ is a graded $\Lambda$-module. Thus its
endomorphism algebra $\End_\Lambda(T_\ww)$ is naturally graded. 
We have the following connection with the previous grading.

\begin{prop}
The isomorphism of algebras $\End_\Lambda(T_\ww)\simeq \Jac(Q_\ww,W_\ww,F)$ of Theorem~\ref{birs} is an isomorphism of graded algebras
 $$(\End_\Lambda(T_\ww),deg)\simeq (\Jac(Q_\ww,W_\ww,F),\varphi)$$ where $deg$ is induced by the grading $deg$ on the preprojective algebra $\Lambda$, and $\varphi$ is the grading on $Q_\ww$ defined in the previous section.
\end{prop}

\begin{proof} 
Each $a:i\rightarrow j$ in $Q_1$ gives maps
$e_i(\Lambda/\Ii_{\ww_r})\rightarrow e_j(\Lambda/\Ii_{\ww_s})$, where $u_r=i$ and $u_s=j$. These maps are obviously of degree $0$
since they are induced by the degree zero map $a:e_i\Lambda\rightarrow
e_j\Lambda$.

Each $a:i\rightarrow j$ in $Q$ induces maps
$e_j(\Lambda/\Ii_{\ww_t})\rightarrow e_i(\Lambda/\Ii_{\ww_s})$, where $u_t=j$ and $u_s=i$. They are induced by the
degree 1 map
$a^*:e_j\Lambda\rightarrow e_i\Lambda$, thus they are maps of
degree $1$.

For any $i$ in $Q_0$, there are surjective maps
$e_i(\Lambda/\Ii_{\ww_t})\rightarrow e_i(\Lambda/\Ii_{\ww_r})$, where $u_t=u_r=i$ and $t>r$. They are induced by the
identity $e_i\Lambda\rightarrow e_i\Lambda$, thus they are maps of
degree 0.

Hence we get the grading $\varphi$ defined in the previous section.
\end{proof} 

\begin{rema} Note that the summands of $\Lambda_w$ are all graded $\Lambda$-modules,
but this does not imply that all the objects in $\Ee_w$ are
gradable. 
In the proof of \cite[Proposition 5.2]{GS}, Geiss and Schr\"oer describe explicitly a non gradable module over the preprojective algebra associated to the Dynkin 
graph $A_6$.\end{rema}

\section{Examples}

In this section we illustrate the previous theory through two examples. The first one is an example covered by 
Theorem~\ref{birs}. It is given by a standard cluster-tilting object in the category $\Ee_w$ for some reduced
 word $\ww$. The second example shows that Theorem~\ref{theoremgrading}  may also apply for cluster-tilting objects in $\Ee_w$
 which are not standard.

\subsection{Standard cluster-tilting object associated to a reduced word}

Let $Q$ be the following graph $$
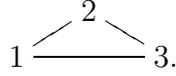
\xymatrix@R=.1cm@C=.5cm{&2\ar@{-}[dr] &
  \\1\ar@{-}[rr]\ar@{-}[ur] &&3.}$$ 
Let $\ww$ be the reduced word $s_1s_2s_3s_1s_3s_2s_1$
in the Coxeter group $C_{Q}$.
An admissible orientation of $Q$ as defined in section 4.2 is $$
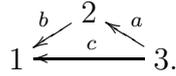
\xymatrix@R=.1cm@C=.5cm{&2\ar@{<-}[dr]^a &
  \\1\ar@{<-}[rr]^c\ar@{<-}[ur]^b &&3.}$$

Let us put the following grading for the preprojective algebra $\Lambda$.

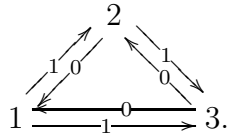
 $$\xymatrix{&2 \ar@<.6ex>[dr]|(.45)1 \ar@<.6ex>[dl]|(.45)0 &
  \\1\ar@<-.6ex>[rr]|(.45)1 \ar@<.6ex>[ur]|(.45)1 &&3.\ar@<-.6ex>[ll]|(.45)0\ar@<.6ex>[ul]|(.45)0}$$
Then the canonical cluster-tilting object $T_\ww$ of the Frobenius
category $\Ee_w$ has the following
indecomposable summands:

\begin{eqnarray*}
T_1={\bsm 1\esm},\quad T_2={\bsm 2\\ 1\esm},\quad T_3={\bsm &3&&\\1&&2&\\
  &&&1\esm},\quad T_4={\bsm &1&&&\\2&&3&&\\&1&&2&\\&&&&1\esm},\\
T_5={\bsm&&3&&\\&1&&2&\\2&&&&1\esm},\quad T_6={\bsm
  &&&2&&&&\\&&3&&1&&&\\&1&&2&&3&&\\2&&&&1&&2&\\&&&&&&&1\esm}\quad
\textrm{and}\quad T_7={\bsm &&&&1&&&\\&&&2&&3&&\\&&3&&1&&2&\\&1&&2&&&&1\\2&&&&&&&\esm}
\end{eqnarray*}
The indecomposable projective injectives are $T_5$, $T_6$ and $T_7$.  
As $\Lambda_0$-module (=$kQ$-module) $T_6$ is isomorphic to the direct sum ${\bsm 2\\3\esm}\oplus {\bsm 1\\2\esm}\oplus{\bsm &1&\\2&&3\esm}\oplus{\bsm 1\esm}\oplus{\bsm 2\esm}\oplus{\bsm 1\esm}$.   

By \cite{Bua2} and \cite{BIRSm}, we know the shape of the quiver of
$B=\End_{\Ee_w}(T)$. Its grading coming from the grading of $\Lambda$ is the following.
$$
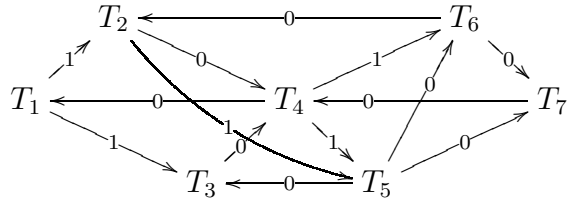
\xymatrix@R=.5cm@C=.5cm{&T_2\ar[drr]|0\ar@/_1pc/[ddrrr]|1&&&&T_6\ar[llll]|0\ar[dr]|0&\\
T_1\ar[ur]|1\ar[drr]|1 &&& T_4\ar[lll]|0\ar[urr]|1\ar[dr]|1&&&T_7\ar[lll]|(.7)0\\
 &&T_3\ar[ur]|(.45)0&&T_5\ar[ll]|0\ar[uur]|(.6)0\ar[urr]|0&&}$$

The algebra $A=B_0$ is then given by the quiver with relations.

$$
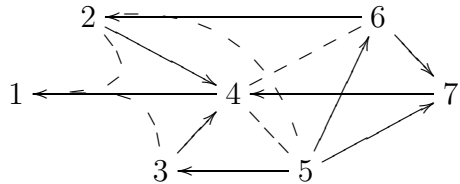
\xymatrix@R=.5cm@C=.5cm{&2\ar[drr]\ar@/^2pc/@{--}[ddrrr]&&&&6\ar[llll]\ar[dr]&\\1\ar@{--}@(r,dr)[ur]\ar@{--}@(r,u)[drr]
  &&& 4\ar[lll]\ar@{--}[urr]\ar@{--}[dr]&&&7\ar[lll]\\
  &&3\ar[ur]&&5\ar[ll]\ar[uur]\ar[urr]&&}$$

The indecomposable projective $A$-modules are

\begin{eqnarray*} {\bsm 1\\4\\7\esm},\quad {\bsm 2\\6\esm}, \quad {\bsm
    3\\5\esm},\quad {\bsm &&4&&\\3&&7&&2\\&5&&6&\esm}, \quad {\bsm
    5\esm},\quad {\bsm 6\\5\esm},\quad \textrm{and}\quad {\bsm
    &&7&\\&6&&5\\5&&&\esm}\end{eqnarray*}
and the indecomposable injectives are 

\begin{eqnarray*} {\bsm1\esm},\quad {\bsm 4\\2\esm},\quad {\bsm
    4\\3\esm}, \quad {\bsm 1\\4\esm}, \quad {\bsm
    &4&&&&7\\3&&7&&6&\\&&5&&&\esm}, \quad {\bsm &4&\\7&&2\\&6&\esm},
  \quad \textrm{and}\quad {\bsm1\\4\\7\esm}\end{eqnarray*}

The algebra $\A$ is given by the quiver with relations:
$$
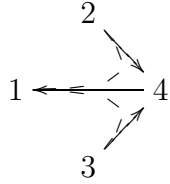
\xymatrix@R=.5cm@C=.5cm{&2\ar[dr]&\\1\ar@{--}@(r,dr)[ur]\ar@{--}@(r,ur)[dr]
  && 4\ar[ll]\\
  &3\ar[ur]&}$$
It is an algebra of global dimension 2. We have 
$\End_{\Cc_{\A}}(\A)\simeq \Jac(\bar{Q}_\ww,\bar{W}_\ww)$ where 

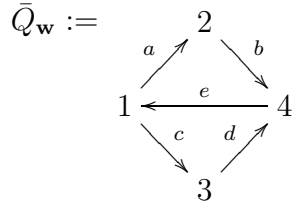
$$\bar{Q}_\ww:=\xymatrix@R=.6cm@C=.6cm{&2\ar[dr]^b&\\1\ar[ur]^a\ar[dr]^c
  && 4\ar[ll]_e\\
  &3\ar[ur]^d&}$$ and $\bar{W}_\ww:=bae+dce$. It
is isomorphic to the algebra $\bar{B}$.

We denote by $G$ the composition 
$$\xymatrix{G:\Dd^{\rm b}(\A)\ar[rr]^{Res} &&\Dd^{\rm b}(A)\ar[rr]^{-\lten_A B} && \Dd^{\rm b}(
  B)\ar[rr]^{-\lten_B T} && \Dd^{\rm b}(\Ee_w)\ar[r] & \Dd^{\rm b}(\Ee_w)/\Dd^{\rm b}(\Pp) \simeq \underline{\Ee}_w.}$$
Let $S_2$ be the simple $\A$-module associated to the vertex $2$. We will show that $G(S_2)$ and $G\circ \mathbb{S}^{-1}(S_2)[2]$ are isomorphic as objects in $\underline{\Ee_w}$, where $\mathbb{S}$ denotes the Serre functor of $\Dd^{\rm b}\A$.

The
restriction of $S_2$ in the category $\Dd^{\rm b}(A)$ is
quasi-isomorphic to the complex
$$\xymatrix{\cdots\ar[r] &0\ar[r] & e_5A\ar[r] &e_6 A\ar[r] & e_2
    A\ar[r] & 0\ar[r] &\cdots}$$
Tensoring with $B$ over $A$ we get the complex
$$\xymatrix{\cdots\ar[r] &0\ar[r] & e_5B\ar[r] &e_6 B\ar[r] & e_2
    B\ar[r] &0\ar[r]& \cdots}.$$
Therefore $G(S_2)$ is the complex
$$\xymatrix{G(S_2)=(T_5\ar[r] &T_6 \ar[r] & T_2)}.$$
Now the simple $\A$-module $S_2$ is quasi-isomorphic in $\Dd^{\rm b}(\A)$ to the complex 
$$\xymatrix{ \cdots \ar[r] & 0\ar[r] & e_2 (D\A)\ar[r] & e_4(D\A)
  \ar[r] & e_1(D\A)\ar[r] & 0\ar[r]&\cdots}$$
Hence the restriction of $S_2\lten_{\A}\RHom_{\A}(D\A,\A)[2]=\mathbb{S}^{-1}(S_2)[2]$ in $\Dd^{\rm b}A$
is isomorphic to 
$$\xymatrix{ \cdots \ar[r] & 0\ar[r] & {\bsm 2\esm}\ar[r] & {\bsm &4&\\2&&3\esm}
  \ar[r] & {\bsm 1\\4\esm}\ar[r] & 0\ar[r]&\cdots}$$
which is quasi-isomorphic to the complex 
$$\xymatrix{ \cdots \ar[r] & 0\ar[r] & {\bsm 5\esm}\oplus{\bsm 2\\6\esm}\ar[r] & {\bsm &&4&&\\3&&7&&2\\&5&&6&\esm}
  \ar[r] & {\bsm 1\\4\\7 \esm}\ar[r] & 0\ar[r]&\cdots}$$
 that is, to the complex
$$\xymatrix{ \cdots \ar[r] & 0\ar[r] & e_5A\oplus e_2 A\ar[r] & e_4 A
  \ar[r] & e_1 A\ar[r] & 0\ar[r]&\cdots}$$
Therefore the object $G(\mathbb{S}^{-1}S_2[2])$ is the complex
$$\xymatrix{G(\mathbb{S}^{-1}S_2[2])=(T_5\oplus T_2\ar[r] &T_4 \ar[r] & T_1)}.$$
Now the simple $\A$-module $S_2$ is also quasi-isomorphic to the complex 
$$\xymatrix{\cdots\ar[r] &0\ar[r] & e_2(DA)\ar[r] &e_4 (DA)\ar[r] &
  e_1 (DA)\ar[r] &0\ar[r] & \cdots}$$
Hence the complex $S_2\lten_{\A} \RHom_A(DA,\A)\lten_A B[2]\simeq \RHom_A(DA,S_2)\lten_A B[2]$ is the complex
$$\xymatrix{\cdots\ar[r] &0\ar[r] & e_2B\ar[r] &e_4 B\ar[r] &
  e_1 B\ar[r] & \cdots}.$$
We have a morphism in the category $\Dd^{\rm b}(B)$:
$$\xymatrix{\cdots\ar[r] &0\ar[r]\ar[d] & e_2B\ar[r]\ar[d] &e_4 B\ar[r]\ar[d] &
  e_1 B\ar[r]\ar[d] & \cdots\\ \cdots\ar[r] &0\ar[r] & e_5B\ar[r] &e_6 B\ar[r] & e_2
    B\ar[r] & \cdots}$$
whose cone is 
$$\xymatrix{\cdots\ar[r] & 0\ar[r] &e_2B\ar[r] & e_5B\oplus e_4B\ar[r] & e_6B\oplus
  e_1B\ar[r]& e_2B\ar[r] & 0\ar[r] & \cdots}.$$ An easy computation shows that it is quasi-isomorphic to the simple
$B$-module $S_2$.
Hence we get the following triangle in $\Dd^{\rm b}(B)$:
$$\xymatrix{S_2\lten_{\A} \RHom_A(DA,\A)\lten_A B[2]\ar[r] & S_2\lten_A B\ar[r]
&S_2\ar[r] &S_2\lten_A \RHom_A(DA,\A)\lten_{\A} B[3]}.$$ 

Note that this triangle is nothing but  $S_2\lten_{\A}(*)$ where $(*)$ is the triangle defined in Proposition~\ref{triangle}.
The object $S_2\lten_BT$ in $\Dd^{\rm b}(\Ee_w)$ is then the complex
$$\xymatrix{\cdots\ar[r] & 0\ar[r] &T_2\ar[r] & T_5\oplus T_4\ar[r] & T_6\oplus
  T_1\ar[r]& T_2\ar[r] & 0\ar[r] & \cdots}.$$
A direct computation shows that it is acyclic. Indeed it is the 2-almost split sequence
associated with $T_2$.
Finally we have morphisms
$$\xymatrix@-.5cm{ & (T_2\rightarrow T_4\rightarrow T_1)\ar[dr]^{(ii)}\ar[dl]_{(i)} & \\
G(S_2)=(T_5\rightarrow T_6\rightarrow T_2)\ar[d]  && (T_5\oplus T_2\rightarrow T_4\rightarrow T_1)=G(\mathbb{S}^{-1}S_2[2])\ar[d]\\ 0\simeq(T_2\rightarrow T_5\oplus T_4\rightarrow T_1\oplus T_6\rightarrow T_2) && (T_5)\in\Dd^{\rm b}(\Pp)}$$
 The cone of the morphism $(i)$ is acyclic, and the cone of the morphism $(ii)$ is $T_5[-2]$, which is perfect. 
Thus in $\underline{\Ee_w}=\Dd^{\rm b}(\Ee_w)/\Dd^{\rm b}(\Pp)$, the objects $G(S_2)$
and $G(\mathbb{S}^{-1}S_2[2])$ are isomorphic.

\subsection{Example which is not associated with a word}

Let $\Ee$ be the category $\mod \Lambda$ where $\Lambda$ is the
preprojective algebra of type $A_3$, which is one of the cases investigated by Geiss, Leclerc and Schr\"oer in \cite{Gei3}. 
This is a Frobenius category which is stably 2-Calabi-Yau and of the form $\Ee_w$ where $w$ is the element in the Coxeter
group of maximal length. Corresponding to the reduced expression $s_1s_2s_3s_1s_2s_1$ is the 
standard cluster-tilting object:
$$
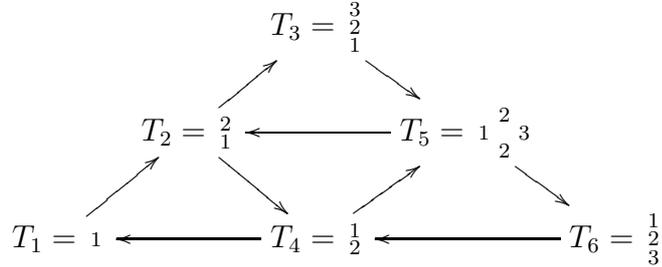
\xymatrix@R=.5cm@C=.2cm{ && T_3={\bsm 3\\2\\1\esm}\ar[dr]&& \\ &T_2={\bsm
    2\\1\esm}\ar[ur]\ar[dr]& &T_5={\bsm&2&\\1&&3\\ &2&\esm}\ar[ll]\ar[dr]
  \\T_1={\bsm 1\esm}\ar[ur] && T_4={\bsm 1\\2\esm}\ar[ll]\ar[ur] &&
  T_6={\bsm1\\2\\3\esm}\ar[ll]}$$
The endomorphism algebra is a frozen Jacobian algebra.

Let us mutate  the object $T_2={\bsm 2\\1\esm}$. Its complement
is $T_2^*={\bsm 3&&1\\&2&\esm}$. The new
cluster-tilting object $T^*$ is given by

$$
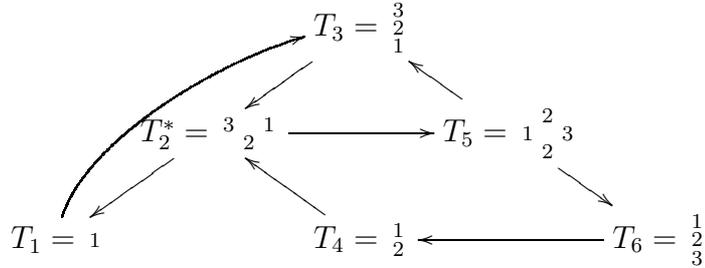
\xymatrix@R=.5cm@C=.2cm{ && T_3={\bsm 3\\2\\1\esm}\ar[dl]&& \\ &T_2^*={\bsm
    3&&1\\&2&\esm}\ar[dl]\ar[rr]& &T_5={\bsm&2&\\1&&3\\ &2&\esm}\ar[ul]\ar[dr]
  \\T_1={\bsm 1\esm}\ar@(u,l)[uurr] && T_4={\bsm 1\\2\esm}\ar[ul] &&
  T_6={\bsm1\\2\\3\esm}\ar[ll]}$$

One can easily check that the endomorphism algebra is isomorphic to the frozen Jacobian
algebra $B=\Jac(Q,W,F)$, where
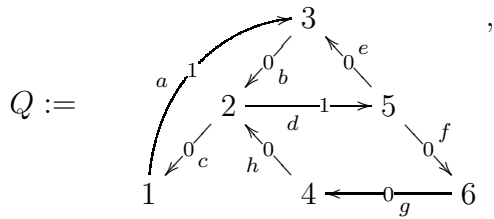
$$\xymatrix@R=.6cm@C=.6cm{ &&& 3\ar[dl]^b|0&& \\Q:=& &2\ar[dl]^c|0\ar[rr]_(.4)d|(.6)1& &5\ar[ul]_e|0\ar[dr]^f|0
  \\&1\ar@(u,l)[uurr]^(.4)a|1 && 4\ar[ul]^h|0 &&
  6\ar[ll]^(.4)g|0},$$ $W:=acb+dbe+dhgf$, $F_0:=\{3,5,6\}$ and $F_1:=\{e,f\}$.

If we put $\varphi(a)=\varphi(d)=1$ and
$\varphi(b)=\varphi(c)=\varphi(e)=\varphi(f)=\varphi(h)=0$, we obtain
a grading satisfying hypothesis (H1)-(H4) of Theorem
\ref{theoremgrading}.
The algebra $A$ is then given by the quiver 
$$
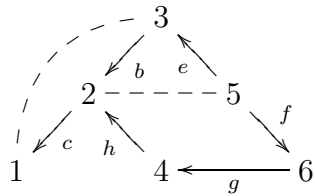
\xymatrix@R=.5cm@C=.5cm{ && 3\ar[dl]^b&& \\ &2\ar[dl]^c\ar@{--}[rr]& &5\ar[ul]^e\ar[dr]^f
  \\1\ar@(u,l)@{--}[uurr] && 4\ar[ul]^h &&
  6\ar[ll]^g}$$ 
with relations $cb=0$ and $be=hgf$.
The algebra $\A$ is the hereditary algebra with quiver
$$\xymatrix{1 & 2\ar[l]^c & 4\ar[l]^h}.$$
One can check that the image of $S_2$ under the functor
$\xymatrix@-.5cm{G:\Dd^{\rm b}(\A)\ar[rr]^{Res} && \Dd^{\rm b}(A)\ar[rr]^{-\lten_A T^*}&&\Dd^{\rm b}(\Ee)}$
is $$\xymatrix{(T_5\ar[rr]^{\left({\bsm e\\fg\esm}\right)} && T_4\oplus
  T_3\ar[rr]^{\left({\bsm b&h\esm}\right)}& & T^*_2).}$$ 
The object $\mathbb{S}^{-1}(S_2)[2]$ is quasi-isomorphic to the complex
$\xymatrix{{\bsm 2\\4\esm}\ar[r] & {\bsm 1\\2\\4\esm}}$, thus its
restriction in $\Dd^{\rm b}(A)$ is quasi-isomorphic to the complex
$$\xymatrix{\cdots \ar[r] &0\ar[r] & e_3 A\ar[r]^{b}
 & e_2 A\ar[r]^c & e_1 A\ar[r] &0\ar[r] & \cdots}.$$ Hence the complex
$G(\mathbb{S}^{-1}S_2[2])$ is
$$\xymatrix{\cdots \ar[r] &0\ar[r] & T_3 \ar[r]^{b}
 & T^*_2 \ar[r]^c & T_1 \ar[r] &0\ar[r] &\cdots}.$$
We have morphisms in $\Dd^{\rm b}(\Ee)$
$$\xymatrix{& (T_2^*\rightarrow T_1)\ar[dl]_{(i)}\ar[dr]^{(ii)} & \\
 G(S_2)=(T_5\rightarrow T_4\oplus T_3\rightarrow T_2^*)\ar[d] && (T_3\rightarrow
 T^*_2\rightarrow T_1)=G\circ\mathbb{S}(S_2)[2]\ar[d]\\ 0\simeq (T_2^*\rightarrow T_1\oplus T_5\rightarrow T_3\oplus T_4\rightarrow T^*_2) & & (T_3)\in\Dd^{\rm b}(\Pp)}$$

The cone of the morphism $(i)$ is $$\xymatrix{ T_2^*\ar[r]^{\left({\bsm
        c\\d\esm}\right)} & T_1\oplus T_5\ar[r]^{\left({\bsm a &e\\0
        &fg\esm}\right)} & T_3\oplus T_4\ar[r]^{\left({\bsm
        b&h\esm}\right)} &T^*_2}$$
which is the 2-almost-split sequence associated to $T^*_2$, hence an acyclic complex.
The cone of the morphism $(ii)$ is $T_3[-2]$, which is perfect.
 Thus in $\Dd^{\rm b}(\Ee)/\Dd^{\rm b}(\Pp)$ the objects $G(S_2)$ and
 $G(\mathbb{S}^{-1}S_2[2])$ are isomorphic.

This example gives some hope that Theorem~\ref{theoremgrading} can be applied to stably 2-CY categories other than those coming from an element of a Coxeter group.

\newcommand{\etalchar}[1]{$^{#1}$}
\providecommand{\bysame}{\leavevmode\hbox to3em{\hrulefill}\thinspace}
\providecommand{\MR}{\relax\ifhmode\unskip\space\fi MR }
\providecommand{\MRhref}[2]{%
  \href{http://www.ams.org/mathscinet-getitem?mr=#1}{#2}
}
\providecommand{\href}[2]{#2}


\end{document}